\documentclass{article}
\usepackage{amsthm}
\usepackage{amsmath}
\usepackage{amssymb}

\title{On Invariant volumes of codimension-one Anosov flows
 and the Verjovsky conjecture}
\author{Masayuki ASAOKA
\footnote{Partially supported by JSPS PostDoctoral Fellowships for
 Research Abroad}\\
Department of Mathematics, Kyoto University\\
 and UMPA, ENS de Lyon\\
}

\def\RR{\mathbb{R}}
\def\ZZ{\mathbb{Z}}

\def\cO{{\cal O}}
\def\cF{{\cal F}}
\def\cG{{\cal G}}
\def\cM{{\cal M}}
\def\cU{{\cal U}}
\def\cR{{\cal R}}
\def\cL{{\cal L}}

\def\--{{\setminus}}

\def\ra{{\rightarrow}}
\def\st{{\;|\;}}

\DeclareMathOperator{\Int}{Int}

\DeclareMathOperator{\Div}{div}
\DeclareMathOperator{\Grad}{grad}

\DeclareMathOperator{\Id}{Id}

\def\vphi{\varphi}
\def\del{\partial}

\def\hsp{{\hspace{3mm}}}

\newcommand\cl[1]{{\overline{#1}}}

\newcommand\ch[1]{{\check{#1}}}

\theoremstyle{plain}
\newtheorem{thm}{Theorem}[section]
\newtheorem{prop}[thm]{Proposition}
\newtheorem{lemma}[thm]{Lemma}
\newtheorem{cor}[thm]{Corollary}
\newtheorem{question}[thm]{Question}
\newtheorem{conj}[thm]{Conjecture}

\newtheorem*{Mthm}{Main Theorem}

\newtheorem*{Ver-conj}{The Verjovsky Conjecture}
\theoremstyle{definition}

\theoremstyle{remark}
\newtheorem{rmk}[thm]{Remark}

\numberwithin{equation}{section}

\begin{document}
\maketitle
\begin{abstract}
We show that any topologically transitive codimension-one Anosov flow
 on a closed manifold is topologically equivalent to
 a smooth Anosov flow that preserves a smooth volume.
By a classical theorem due to Verjovsky,
 any higher dimensional codimension-one Anosov flow
 is topologically transitive.
Recently, Simi\'c showed that
 any higher dimensional codimension-one Anosov flow
 that preserves a smooth volume
 is topologically equivalent to the suspension
 of an Anosov diffeomorphism.
Therefore, our result gives a complete classification
 of codimension-one Anosov flow up to topological equivalence
 in higher dimensions.
\end{abstract}

\section{Introduction}
%
%
In \cite{LS},
 Liv\v sic and {Sina\u \i } showed that
 an Anosov flow $\Phi=\{\Phi^t\}_{t \in \RR}$ on a closed manifold
 preserves an absolutely continuous measure
 with respect to the Lebesgue measure
 if and only if $\det (D\Phi^T)_p=1$
 for any periodic point $p$ with period $T$.
In \cite{LMM},
 de la Llave, Marco, and Moriy\'on showed
 the corresponding result for smooth volumes.
Their results give a characterization of Anosov flows
 with $C^\infty$ invariant volume.
However, the following natural question is still open and important.
\begin{question}
Which Anosov flow is topologically equivalent to
 a volume preserving one?
\end{question}
We answer this question for codimension Anosov flows. That is,
\begin{Mthm}
Any topologically transitive codimension-one Anosov flow
 on a closed manifold is topologically equivalent to
 a $C^\infty$ Anosov flow with a $C^\infty$ invariant volume.
\end{Mthm}
It is well-known that any Anosov flow that preserves a volume
 is topologically transitive.
Therefore, the above theorem is optimal for codimension-one case.
Remark that there exists a three-dimensional Anosov flow
 that is not topologically transitive (see \cite{FW}).

When the dimension of the manifold is greater than three,
 the main theorem yields an important conclusion.
In 1970's, Verjovsky conjectured that
 any higher dimensional codimension-one Anosov flow should
 be a classical one. That is,
\begin{Ver-conj}
Any codimension-one Anosov flow
 on a closed manifold of dimension greater than three
 is topologically equivalent to the suspension flow
 of a toral automorphism.
\footnote{Ghys \cite{Gh}
 was the first literature that mentioned the conjecture.
 But, he has pointed out that Verjovsky had proposed 
 the conjecture in 1970's.}
\end{Ver-conj}
We recall two results on codimension-one Anosov flows
 on higher dimensional manifolds.
Let $M$ be a closed manifold of dimension greater than three
 and $\Phi$ be a $C^\infty$ codimension-one Anosov flow on $M$.
\begin{thm}
[Verjovsky \cite{Ve}, see also \cite{Ba}]
 $\Phi$ is topologically transitive.
\end{thm}
\begin{thm}
[Simi\'c \cite{Si}] 
If $\Phi$ preserves a $C^\infty$ volume,
 then it is topologically equivalent
 to the suspension flow of a hyperbolic toral automorphism.
\end{thm}
With their results, the main theorem implies
\begin{cor}
The Verjovsky conjecture is true. 
\end{cor}

As far as the author's knowledge,
 no counterexample is known to the analogy of
 the main theorem for higher codimension.
So, we pose
\begin{conj}
 Any topologically transitive Anosov flow is topologically
 equivalent to a volume preserving one.
\end{conj}

The main ingredient of the proof of the main theorem is
  a generalization of Cawley's deformation of Anosov systems.
In \cite{Ca}, Cawley gave a method to deform
 a two-dimensional Anosov diffeomorphism
 into another one that has the desired expansions
 along the stable foliation and the unstable foliation.
In Section \ref{section:Radon},
 we generalize the Radon-Nikodym realization theorem,
 which played the central role in her method,
 to codimension-one Anosov flows.
It allows us to deform an Anosov flow
 into another one that preserves a volume in some sense,
 but the deformation destroys the smoothness of the flow.
So, we need more effort to obtain a smooth flow.
It is done in Section \ref{section:main thm}.

\paragraph{Acknowledgments}
This paper was written while the author stayed at
 Unit\'e de Math\'ematiques Pures et Appliqu\'ees,
 \'Ecole Normale Sup\'erieure de Lyon.
He would like to thank the members of UMPA, especially Professor
 \'Etienne Ghys, for their warm hospitality.
The author would also like to thank Shigenori Matsumoto,
 who gave many suggestions to improve the paper.
Especially, he pointed out that Theorem \ref{thm:foliation}
 is due to Hart and he gave a simpler proof of Lemma \ref{lemma:div}.

%
%
%
\section{Preliminaries}

%
%
\subsection{Notations}
By $\ZZ$ and $\RR$, we denote the set of integers and real numbers
 respectively.
For an integer $r \geq 0$ and $\alpha \in (0,1)$,
 we say a map between $C^\infty$ manifolds is of class $C^{r+\alpha}$
 if it is of class $C^r$ and
 all $r$-th partial derivatives are $\alpha$-H\"older continuous.
We say a map is of class $C^{r+}$ if
 it is of class $C^{r+\alpha}$ for some $\alpha \in (0,1)$.

Let $\pi:E \ra B$ be a finite-dimensional vector bundle
 over a Hausdorff space $B$
 and $g$ be a continuous metric on $E$.
We denote the norm corresponding to $g$ by $\|\cdot\|_g$
 and the fiber $\pi^{-1}(p)$ of $p \in B$ by $E(p)$.
For a linear map $F$ between fibers $E(p)$ and $E(q)$,
 we denote the determinant of $F$
 with respect to the metric $g$ by $\det_g F$.
We define {\it the operator norm} $\|F\|_g$ of $F$ by
\begin{equation*}
\|F\|_g=\sup\left\{
 \|F(v)\|_g \st v \in E(p), \|v\|_g \leq 1 \right\}.
\end{equation*}

\subsection{Regularity of foliations}
Let $\cF$ be a continuous foliation on a manifold $M$.
We denote the leaf that contains $p \in M$ by $\cF(p)$.
For an open subset $U$ of $M$,
 let $\cF|_U$ be the foliation on $U$ such that
 $(\cF|_U)(p)$ is the connected component of $\cF(p) \cap U$
 containing $p \in M$.
A coordinate $\vphi=(x_1,\cdots, x_n)$ on $U$
 is called {\it a foliation coordinate} of $\cF$
 if $x_{m+1},\cdots,x_n$ are constant functions
 on each leaf of $\cF|_U$, where $m$ is the dimension of $\cF$.
A foliation is of class $C^{r+}$ if it is covered by
 $C^{r+}$ foliation coordinates.
We denote the tangent bundle of $M$ by $TM$.
If $\cF$ is a $C^1$ foliation, then we denote the tangent bundle
 of $\cF$ by $T\cF$.
For a homeomorphism $h$ from $M$ to another manifold $M'$,
 we define a foliation $H(\cF)$ on $M'$
 by $H(\cF)(p)=H(\cF(H^{-1}(p)))$.

In general, if a foliation $\cF$ is of class $C^{r+}$,
 then $T\cF$ is of class $C^{(r-1)+}$.
However, the following theorem due to Hart implies that
 we may replace $\cF$ in its $C^{r+}$-equivalence class
 so that $T\cF$ is of class $C^{r+}$.
\begin{thm}[Hart {\cite[Theorem B]{Ha}}]
\label{thm:foliation}
For any $C^{r+\alpha}$-foliation $\cF$
 on a $C^\infty$ closed manifold $M$,
 there exists a $C^{r+\alpha}$ diffeomorphism $H$ of $M$
 such that $T(H(\cF))=DH(T\cF)$ is a $C^{r+\alpha}$ subbundle of $TM$.
Moreover, we can choose the diffeomorphism $H$
 so that it is arbitrary $C^r$-close to the identity map.
\end{thm}
We give a short proof in Appendix \ref{sec:A foliation}.
The theorem has an immediate corollary, which is important for our proof
 of the main theorem.
\begin{cor}
\label{cor:foliation} 
Let $M_1$ be a $C^{1+}$ closed manifold
 and $\cL$ an oriented $C^{1+}$ one-dimensional foliation on $M_1$.
Then, $M_1$ admits a $C^{\infty}$-structure 
 that is compatible with the $C^{1+}$-structure of $M_1$
 and such that $\cL$ is generated by a $C^{1+}$ vector field.
\qed
\end{cor}

%
%
\subsection{$C^\infty$ invariant volumes}
We show the following proposition in this subsection.
\begin{prop}
\label{prop:volume}
Let $X_0$ be a $C^{1+}$ vector fields on a $C^\infty$ 
 closed manifold $M$.
If the flow generated by $X_0$ preserves a H\"older continuous volume,
 then $X_0$ can be $C^1$-approximated by
 a $C^\infty$ vector field that generates
 a flow with a $C^\infty$ invariant volume.
\end{prop}

For a $C^1$ vector field $X$ on $M$
 and a continuous metric $g$ on $TM$,
 we define a function $\Div_g X$ on $M$ by
\begin{equation*}
 \Div_g X(p)=\lim_{t \ra 0}\frac{1}{t}
 \log {\det}_g D\Phi^t(p)
\end{equation*}
 if the limit exists for any $p \in M$,
 where $\Phi$ is the $C^1$ flow generated by $X$.
Remark that the formula
\begin{equation}
\label{eqn:Div}
 {\Div}_{e^h \cdot g} X=(\dim M) \cdot Xh +{\Div}_g X
\end{equation}
 holds when $\Div X$ is well-defined and $h$ is differentiable
 along $X$.
\begin{lemma}
\label{lemma:div}
Let $n$ be the dimension of $M$.
For any $C^{1+\alpha}$ vector field $X$
 and any $C^\infty$ metric $g$ on $TM$,
 the function $\Div_g X$ is well-defined and of class $C^\alpha$.
Moreover,
 if a sequence $\{X_k\}_{k \geq 1}$ of $C^{1+\alpha}$ vector fields
 converges to $X$ with respect to the $C^{1+\alpha}$-topology,
 then $\{\Div_g X_k\}_{k \geq 1}$ converges to $\Div_g X$
 with respect to the $C^\alpha$-topology.
\end{lemma}
\begin{proof}
The metric $g$ induces a $C^\infty$ volume form $\omega_g$ on $M$.
Fix $p \in M$.
By Morser's lemma,
 there exists a $C^\infty$ coordinate $\vphi=(x_1,\cdots,x_n)$
 on an open neighborhood $U$ of $p$ such that
 $(\vphi^{-1})^*\omega_g$ is the standard volume form
 $dx_1 \wedge \cdots \wedge d x_n$.
We define functions $a_1,\cdots,a_n$ on $\vphi(U)$ by
\begin{equation*}
 D\vphi(X)=\sum_{i=1}^n a_n\frac{\del}{\del x_i}.
\end{equation*}
The standard argument gives a formula
\begin{equation*}
\Div_g X=\sum_{i=1}^n \frac{\del a_i}{\del x_i}
\end{equation*}
 even if $X$ is only of class $C^{1+\alpha}$.
Hence, $\Div_g X$ is a well-defined $C^\alpha$ function
 and it depends continuously on $X$.
\end{proof}

Now, we prove Proposition \ref{prop:volume}.
Let $X_0$ be a $C^{1+\alpha}$ vector field on $M$
 that preserves a H\"older continuous volume and
 $\Phi_0$ be the $C^{1+\alpha}$ flow generated by $X_0$.
Fix a $C^\infty$ metric $g_0$ on $TM$.
There exists a $C^\alpha$ function $h_0$
 such that $\det_{e^{h_0} \cdot g_0} D\Phi_0^t(p)=1$
 for any $p \in M$ and $t \in \RR$.
By Lemma \ref{lemma:div},
 $\Div_{g_0} X$ is well-defined and of class $C^\alpha$.
Since
\begin{equation*}
 (\dim M) \cdot (h_0 \circ \Phi_0^t(p)-h_0(p))
 +\log {\det}_{g_0} D\Phi_0^t(p)
 =\log {\det}_{e^{h_0} \cdot g_0} D\Phi_0^t(p)=0,
\end{equation*}
$h_0$ is differentiable along $X_0$ and
\begin{equation}
\label{eqn:volume-1}
(\dim M) \cdot X_0 h_0=-\Div_{g_0} X_0. 
\end{equation}
The standard construction by using a mollifier gives
 a sequence $\{h_k\}_{k \geq 1}$ of $C^\infty$ functions on $M$
 such that
\begin{equation*}
\|h_k-h_0\|_{C^\alpha}+ \|X_0h_k-X_0h_0\|_{C^\alpha} <\frac{1}{k},
\end{equation*}
 where $\|\cdot\|_{C^\alpha}$ is the $C^\alpha$-norm of a function.
Take a sequence $\{Y_k\}_{k \geq 1}$ of
 $C^\infty$ vector fields on $M$ that converges
 to $X_0$ with respect to the $C^{1+\alpha}$-topology.
By Lemma \ref{lemma:div},
 for any fixed $k$, $\Div_{e^{h_k}g_0}Y_{k'}$ converges to
 $\Div_{e^{h_k}g_0}X_0$ as a $C^\alpha$ function
 when $k'$ tends to infinity.
Hence, by taking a subsequence of $\{Y_k\}_{k \geq 1}$ if it is necessary,
  we may assume that
\begin{equation*}
 \|\Div_{e^{h_k} \cdot g_0} Y_k - \Div_{e^{h_k} \cdot g_0} X_0\|_{C^\alpha}
 <\frac{1}{k} 
\end{equation*}
By the equations (\ref{eqn:Div}) and (\ref{eqn:volume-1}),
 we have
\begin{equation*}
\Div_{e^{h_k} \cdot g_0} Y_k
 = \left(\Div_{e^{h_k} \cdot g_0} Y_k
    -\Div_{e^{h_k} \cdot g_0} X_0\right)
 + (\dim M) \cdot \left(X_0 h_k-X_0 h_0 \right)
\end{equation*}
  and the right-hand term converges to $0$.

Now, we mimic the proof of Theorem 2.2 in \cite{Ar}
 due to Arbeito and Matheus.
We fix a point $p_* \in M$.
Let $\Delta_k$ and $\Grad_k f$
 be the Laplacian and the gradient flow of a $C^\infty$ function $f$
 with respect to the metric $e^{h_k} \cdot g_0$.
Since the integral of $-\Div_{e^{h_k} \cdot g_0}$ over $M$
 is zero,
 the partial differential equation
 $\Delta_k f_k= -\Div_{e^{h_k} \cdot g_0} Y_k$
 has a unique $C^\infty$ solution $f_k$ with $f_k(p_*)=0$
 and the sequence $(f_k)_{k \geq 1}$ converges to $0$
 in the $C^{2+\alpha}$-topology.
The vector field $X_k=Y_k+\Grad_k f_k$
 satisfies
\begin{equation*}
 \Div_{e^{h_k} \cdot g_0} X_k=
 \Div_{e^{h_k} \cdot g_0} Y_k
 + \Delta_k f_k=0
\end{equation*}
 and the sequence $\{X_k\}_{k \geq 1}$
 converges to $X_0$ with respect to the $C^{1+\alpha}$-topology.

\subsection{Anosov flows}
Let $M$ be a $C^\infty$ manifold 
 and $\Phi$ a $C^{1+}$ flow without stationary points.
By $T\Phi$, we denote the one-dimensional subbundle of $TM$
 that is tangent to the orbits of $\Phi$.
We say $\Phi$ is {\it Anosov} if there exists
 a continuous $D\Phi$-invariant splitting
 $TM=T\Phi \oplus E^{ss} \oplus E^{uu}$,
 a continuous metric $g$ on $TM$,
 and a constant $\lambda>0$ such that
\begin{equation}
\label{eqn:Anosov ss,uu}
\max\left\{
 \|D\Phi^t|_{E^{ss}(p)}\|_g, \|D\Phi^{-t}|_{E^{uu}(p)}\|_g
 \right\} < \exp(-\lambda t)
\end{equation}
 for any $p \in M$ and any $t > 0$.
The splitting $TM=T\Phi \oplus E^{ss} \oplus E^{uu}$
 is called {\it an Anosov splitting} of $\Phi$.
It is known that the splitting is H\"older continuous
 and the subbundles
 $E^{ss}$, $E^{uu}$, $T\Phi \oplus E^{ss}$, $T\Phi \oplus E^{uu}$
 are uniquely integrable.
The corresponding foliations are called
 {\it the strong stable foliation},
 {\it the strong unstable foliation},
 {\it the weak stable foliation}, and
 {\it the weak unstable foliation}, respectively.

We say $\Phi$ is a {\it codimension-one} Anosov flow
 if $E^{uu}$ is one-dimensional.
It is known that
 if $\Phi$ is a $C^\infty$ codimension-one Anosov flow,
 then $T\Phi \oplus E^{ss}$ is of class $C^{1+}$.
See {\it e.g.} \cite[Corollary 19.1.12]{KH}.

\subsection{The Gibbs measure for subshift of finite types}
In this subsection,
 we review some basic results on symbolic dynamics.
Fix a positive integer $k_*$.
Let $\Sigma$ be the set of maps from
 $\{m \in \ZZ \st m \geq 0\}$ to $\{1,\cdots,k_*\}$.
We define a distance $d_\Sigma$ on $\Sigma$ by
\begin{equation}
\label{eqn:d-Sigma}
 d_\Sigma(\xi,\xi') =\left\{
\begin{array}{ll}
0 & (\xi=\xi') \\
\exp\left(\inf\{m \geq 0 \st \xi(m) \neq \xi'(m)\}\right)
 & (\xi \neq \xi').
\end{array}
\right.
\end{equation}
We also define {\it the shift map} $\sigma$ by
 $\sigma(\xi)(m)=\xi(m+1)$.
The metric space $(\Sigma,d_\Sigma)$ is compact
 and the map $\sigma$ is Lipschitz continuous.

Fix a $(k_* \times k_*)$-matrix $A=(a_{ij})$
 with entries in $\{0,1\}$.
Suppose that $A$ is {\it mixing},
 that is, there exists $m_0 \geq 1$ such that
 all entries of $A^{m_0}$ are positive.
Let $\Sigma_A$ be the set of maps $\xi \in \Sigma$ such that
 $a_{\xi(k+1)\xi(k)}=1$ for any $k \geq 0$.
It is a $\sigma$-invariant closed subspace of $\Sigma$.
By $\sigma_A$, we denote the restriction of $\sigma$ to $\Sigma_A$.
We call the pair $(\Sigma_A, \sigma_A)$
 {\it the subshift of finite type} associated with the matrix $A$.

Let $\pi_{A,m}$ be the natural projection from
 $\Sigma_A$ to the set of maps
 from $\{0,\cdots,m-1\}$ to $\{1,\cdots,k_*\}$.
Put $\Sigma_{A,m}=\pi_{A,m}(\Sigma_A)$.
For any given H\"older continuous function $g$ on $\Sigma_A$,
 we define {\it the topological pressure} $P_{\sigma_A}(g)$ of $g$
 with respect to $\sigma_A$ by
\begin{equation}
\label{eqn:pressure}
P_{\sigma_A}(g)=
\lim_{m \ra +\infty} \frac{1}{m}\left[
 \log \sum_{\xi_m \in \Sigma_{A,m}}
 \exp\left(\sup_{\xi \in \pi_{A,m}^{-1}(\xi_m)}
 \sum_{j=0}^{m-1} g \circ \sigma_A^j(\xi)\right)\right].
\end{equation}
It satisfies
\begin{align}
\label{eqn:pressure-2-1}
P_{\sigma_A}(g_1)-P_{\sigma_A}(g_2)
 &\leq
 P_{\sigma_A}(0) \cdot \sup_{\xi \in \sigma_A} (g_1(\xi)-g_2(\xi)),\\
\label{eqn:pressure-2-2}
P_{\sigma_A^m}\left(\sum_{j=0}^{m-1}g_1 \circ \sigma_A^j\right)
 & = m \cdot P_{\sigma_A}(g_1)
\end{align}
 for any H\"older continuous functions $g_1$, $g_2$,
 and any integer $m \geq 1$.
Since $A$ is mixing,
 we can see that $P_{\sigma_A}(0)$ is positive.

\begin{lemma}
\label{lemma:pressure}
Let $g$ be a H\"older continuous function on $\Sigma_A$
 such that
\begin{equation*}
 \inf_{\xi \in \Sigma_A} \sum_{j=0}^{m-1} g \circ \sigma_A^j(\xi) >0
\end{equation*}
 for some $m \geq 1$.
Then, there exists $\rho=\rho(g)>0$ such that
 $P_{\sigma_A}(-\rho \cdot g)=0$.
\end{lemma}
\begin{proof}
Put $g_m=\sum_{j=0}^{m-1} g \circ \sigma_A^j$.
By the assumption,
 there exists $C>1$ such that $C^{-1} \leq g_m(\xi) \leq C$
 for any $\xi \in \Sigma_A$.
By (\ref{eqn:pressure-2-2}), we have
 $P_{\sigma_A^m}(-\rho \cdot g_m)
  =m \cdot P_{\sigma_A}(-\rho  \cdot g)$ for any $\rho>0$.
The inequality (\ref{eqn:pressure-2-1}) implies
\begin{equation*}
C^{-1}
 \leq
 \frac{P_{\sigma_A^m}(-\rho_2 \cdot g_m)
  -P_{\sigma_A^m}(-\rho_1 \cdot g_m)
  }{P_{\sigma_A^m}(0) \cdot (\rho_1-\rho_2)}
  \leq C
\end{equation*}
 for any $\rho_1>\rho_2>0$.
Hence, the function $\rho \mapsto P_{\sigma_A}(-\rho \cdot g)$
 is continuous, strictly decreasing,
 and unbounded from the below.
Since $P_{\sigma_A}(0)>0$, 
 there exists $\rho>0$ such that $P_{\sigma_A}(-\rho \cdot g)=0$.
\end{proof}

Let $X,Y$ be topological spaces and $f:X \ra Y$ be a covering map.
For a Borel measure $\mu$ of $Y$,
 there exists a unique Borel measure $\nu$ of $X$
 that satisfies $\nu(A \cap U)=\mu(h(A \cap U))$
 for any Borel subset $A$ of $X$
 and any open subset $U$ such that $h|_{U}$ is
 an homeomorphism onto its image.
We call the measure $\nu$ {\it local push-forward} of $\mu$
 and denote it by $\mu \circ h$.
Notice that $\sigma_A$ is a local homeomorphism,
 and hence, we can define $\mu \circ \sigma_A$ 
 for any Borel measure $\mu$ on $\Sigma_A$.

\begin{thm}
\label{thm:Gibbs}
For any given H\"older function $g$ on $\Sigma_A$,
 there exists a Borel probability measure $\mu_g$
 and a positive H\"older continuous function $h_g$ on $\Sigma_A$
 such that
\begin{enumerate}
 \item
 $\mu_g$ is non-atomic and
 positive on each non-empty open subset of $\Sigma_A$,
 \item
 the measure $h_g \cdot \mu_g$ is $\sigma_A$-invariant and ergodic,
 and
 \item for any $\xi \in \Sigma_A$,
\begin{equation*}
 \log\frac{d (\mu_g \circ \sigma_A)}{d \mu_g}(\xi)
 =-g(\xi)+P_{\sigma_A}(g)
\end{equation*}
\end{enumerate}
\end{thm}
\begin{proof}
See Chapter 1 of \cite{Bo}.
\end{proof}

\section{The Radon-Nikodym Realization Theorem}
\label{section:Radon}

%
%

Fix a $C^\infty$ $n$-dimensional closed manifold $M$
 and a $C^\infty$ codimension-one
 topologically transitive Anosov flow  $\Phi$ on $M$.
Let $TM=T\Phi \oplus E^{ss} \oplus E^{uu}$ be
 the Anosov splitting of $\Phi$.
By $\cF^{ss}$, $\cF^{uu}$, $\cF^s$, and $\cF^u$,
 we denote the strong stable foliation,
 the strong unstable foliation,
 the weak stable foliation, and
 the weak unstable foliation of $\Phi$, respectively.
Fix a $C^\infty$ Riemannian metric $g$ on $M$
 that satisfying the equation (\ref{eqn:Anosov ss,uu})
 for some $\lambda>0$.
We denote the distance induced from $g$ by $d$.
Then, we have
\begin{equation}
\label{eqn:Anosov F-ss}
 \limsup_{t \ra +\infty}\frac{1}{t}
 \log d(\Phi^t(q_1),\Phi^t(q_2)) \leq -\lambda
\end{equation}
 for any $q_1 \in M$ and $q_2 \in \cF^{ss}(q_1)$.

Let $B(p,\epsilon)$ be the open $\epsilon$-ball centered at $p \in M$.
For $p \in M$, $\delta>0$, and $\sigma=s,u,ss,uu$,
 let $\cF_\delta^\sigma(p)$ be the connected component
 of $\cF^\sigma(p) \cap B(p,\delta)$ that contains $p$.
There exists $0<\delta_0<\delta_1<\delta_2$ such that
 $\cF_{\delta_2}^\sigma(p)$ is a disk for any $p \in M$
 and any $\sigma=s,u,ss,uu$, and
 the holonomy map
 $h_{pq}:\cF^{uu}_{\delta_1}(q) \ra \cF^{uu}_{\delta_2}(p)$
 of the foliation $\cF^s$
 is well-defined for any $p \in M$ and $q \in B(p,\delta_0)$.

Let $\cM(\Phi)$ be the set of families $\nu=\{\nu_p\}_{p \in M}$
 of measures that satisfy the following properties:
\begin{enumerate}
\item $\nu_p$ is a Borel measure on $\cF^{uu}(p)$
 that is locally finite, non-atomic, and positive
 on each non-empty open subset of $\cF^{uu}(p)$.
\item If $\cF^{uu}(p)=\cF^{uu}(q)$, then
 $\nu_p=\nu_q$.
\item The Radon-Nikodym derivative $d(\nu_p \circ h_{pq})/d\nu_q(q)$
 at $q$ is well-defined for any $p \in M$ and $q \in B(p,\delta_0)$,
 and it is H\"older continuous with respect to $p$ and $q$.
\end{enumerate}

The rest of this section is devoted to
 the proof of the following theorem, which is
 a keystone of the proof of the main theorem.
\begin{thm}
[The Radon-Nikodym Realization Theorem]
\label{thm:Radon}
For any H\"older continuous positive function $f$ on $M$,
 there exists $\nu=\{\nu_p\}_{p \in M} \in \cM(\Phi)$
 and $\rho>0$ such that
\begin{displaymath}
\log\left(\frac{d(\nu_{\Phi^t(p)} \circ \Phi^t)}{d\nu_p}\right)(p)
 =\rho \cdot \int_0^t f\circ \Phi^\tau(p) d\tau
\end{displaymath}
 for any $p \in M$ and $t \in \RR$.
\end{thm}
\begin{rmk}
The original version of the theorem is shown by Cawley \cite{Ca}
 for two dimensional Anosov diffeomorphisms.
Our proof follows her argument, but we construct a family of
 measures on leaves of $\cF^{uu}$ directly,
 not only a transverse measure class of $\cF^s$.
\end{rmk}

\subsection{Markov partitions}
We call a quadruple $R=(U,\vphi,R^s,R^u)$
 {\it a rectangle} for $\Phi$ if
\begin{itemize}
\item $U$ is an open subset of $M$,
\item $\vphi$ is a continuous coordinate on $U$ satisfying
\begin{align*}
\vphi^{-1}(\RR \times \RR^{n-2} \times y)
 & \subset \cF^s(\vphi^{-1}(w,x,y)),\\
\vphi^{-1}(\RR \times x \times \RR)
 & \subset \cF^u(\vphi^{-1}(w,x,y))
\end{align*}
 for any $(w,x,y) \in \vphi(U) \subset \RR^n
 =\RR \times \RR^{n-2} \times \RR$, and
\item $R^s$ and $R^u$ are compact subsets of $\RR^{n-2}$ and $\RR$
 respectively, and $[0,1] \times R^s \times R^u \subset \vphi(U)$.
\end{itemize}
For a rectangle $R=(U,\vphi,R^s,R^u)$,
 we define subsets $\cl{R}$, $\Int R$, $\del_-' R$, and $\del_+' R$
 of $M$ by
\begin{align*}
 \cl{R} &= \vphi^{-1}([0,1] \times R^s \times R^u),\\
 \Int R &= \vphi^{-1}((0,1) \times \Int R^s \times \Int R^u),\\
 \del_-' R &= \vphi^{-1}(0 \times \Int R^s \times \Int R^u),\\ 
 \del_+' R &= \vphi^{-1}(1 \times \Int R^s \times \Int R^u).
\end{align*}
A rectangle $R$ is called {\it proper}
 if $\cl{\Int R^s}=R^s$ and $\cl{\Int R^u}=R^u$.
We call a finite family
 $\cR=\{R_i=(U_i,\vphi_i,R^s_i,R^u_i)\}_{i=1}^{i_*}$
 of proper rectangles
 {\it a Markov partition} associated with $\Phi$ if
 $M=\bigcup_{i=1}^{i_*} \cl{R}_i$,
 $\Int R_i \cap \Int R_j=\emptyset$ for $i \neq j$, and
\begin{align}
\label{eqn:Markov s}
\vphi_j^{-1}(1 \times R^s_j \times y_j) & \subset
 \vphi_i^{-1}(0 \times R^s_i \times y_i),\\
\label{eqn:Markov u}
\vphi_i^{-1}(0 \times x_i \times R^u_i) & \subset
 \vphi_j^{-1}(1 \times x_j \times R^u_j)
\end{align}
 for any $i,j=1,\cdots,i_*$
 and $p=\vphi_j^{-1}(1,x_j,y_j)=\vphi_i^{-1}(0,x_i,y_i)
 \in \del_+'R_j \cap \del_-' R_i$.
{\it The transition matrix} of $\cR$ is
 a $(i_* \times i_*)$-matrix $(a_{ij})$ that is defined by
 $a_{ij}=1$ if $\del'_- R_i \cap \del'_+ R_j \neq \emptyset$
 and $a_{ij}=0$ otherwise.
The following theorem is well-known (see \cite{Ra}).
\begin{thm}
\label{thm:Markov}
For any $\delta>0$, there exists a Markov partition
 $\cR=\{R_i=(U_i,\vphi_i,R^s_i,R^u_i)\}_{i=1}^{i_*}$ such that
 its transition matrix is mixing
 and the diameter of $U_i$ is
 less than $\delta$ for any $i=1,\cdots,i_*$.
\end{thm}

Fix a Markov partition
 $\cR=\{R_i=(U_i,\vphi_i,R^s_i,R^u_i)\}_{i=1}^{i_*}$
 such that the diameter of $U_i$ is less than $\delta_0/8$
 and $\del'_- R_i \cap \del'_+ R_i=\emptyset$ for any $i=1,\cdots,i_*$.
Remark that the holonomy map
 $h_{pq}:\cF^{uu}_{\delta_1}(q) \ra \cF^{uu}_{\delta_2}(p)$
 is well defined for any $p \in U_i$ and $q \in U_j$
 with $U_i \cap U_j \neq \emptyset$.

Let $\Sigma_A$ be the subshift of finite type
 associated with the transition matrix of $\cR$.
For $p \in M$, a pair $(\tau_p,\xi_p)$ is called {\it an itinerary}
 of $p$ if $\tau_p$ is a strictly increasing function
 on $\{m \in \ZZ \st m \geq 0\}$ with $\tau_p(0)=0$,
 $\xi_p$ is an element of $\Sigma_A$, and
 $\Phi^t(p) \in \cl{R}_{\xi(m)}$
 for any $m \geq 0$ and $t \in [\tau_p(m),\tau_p(m+1)]$.
Any $p \in M$ admits at least one itinerary.

Fix a point $p_i$ of $\Int R_i$ for each $i=1,\cdots, i_*$.
Put $(w_i,x_i,y_i)=\vphi_i(p_i)$
 and $I^u_i=\vphi_i^{-1}(w_i \times x_i \times \RR)$.
We define a map $\pi_i$ from
 $U_i$ to $I^u_i$
 by $\pi_i(\vphi_i^{-1}(w,x,y))=\vphi_i^{-1}(w_i,x_i,y)$.
In other words, $\pi_i$ is the projection to $I^u_i$
 along $\cF^s$.
We put $\Lambda^u_i=\vphi_i^{-1}(w_i \times x_i \times R^u_i)$
 and $\Lambda^u=\bigcup_{i=1}^{i_*} \Lambda^u_i$.

Let $\Lambda^u_{ij}$ be the subset of $\Lambda^u$
 consisting of points $p$ that admits
 an itinerary $(\tau_p,\xi_p)$ with $\xi_p(0)=j$
 and $\xi_p(1)=i$.
We define a map $\Phi_{ij}$ from $\Lambda^u_{ij}$ to $\Lambda^u_i$
 by $\Phi_{ij}(p)=\pi_i \circ \Phi^{\tau_p(1)}(p)$.
By the standard argument,
 we can show that for any given $\xi \in \Sigma_A$,
 there exists a unique $p_\xi \in \Lambda^u$
 that admits an itinerary $(\tau_p,\xi_p)$ with $\xi_p=\xi$.
We define a map $\pi_A:\Sigma_A \ra \Lambda^u$ by
 $\pi_A(\xi)=p_\xi$.
Since any $p \in M$ admits an itinerary,
 the map $\pi_A$ is surjective.
The uniqueness of $p_\xi$ implies that
\begin{equation*}
 \pi_A \circ \sigma_A(\xi)=\Phi_{\xi(1)\xi(0)}(\pi_A(\xi))
\end{equation*}
 for any $\xi \in \Sigma_A$, and hence,
\begin{equation}
\label{eqn:pi-A}
\pi_A \circ \sigma_A^m(\xi_p)=\pi_{\xi_p(m)} \circ \Phi^{\tau_p(m)}(p)
\end{equation}
 for any $p \in M$, any itinerary $(\tau_p,\xi_p)$ of $p$,
 and any $m \geq 0$.

Recall that $\Sigma_A$ is a metric space with the metric $d_\Sigma$
 that is given by (\ref{eqn:d-Sigma}).
The set $\Lambda^u$ admits a natural metric $d_0$
 as a subset of the Riemannian manifold $M$.
\begin{lemma}
\label{lemma:pi-A}
The map $\pi_A$ is H\"older continuous.
\end{lemma}
\begin{proof}
Since $\Phi$ expands the foliation $\cF^{uu}$,
 there exists $\lambda'>0$ and a distance $d'$ on $\Lambda^u$ such that
 $d'$ is comparable with $d_0$ and
 $d'(\Phi_{ij}(p),\Phi_{ij}(q)) \geq e^{\lambda'} d'(p,q)$
 for any $p,q \in \Lambda^u_{ij}$.
Let $C$ be the diameter of $\Lambda^u$ with respect to $d'$.
Since
\begin{equation*}
 \pi_A \circ  \sigma_A^m(\xi) =
 \Phi_{\xi(m)\xi(m-1)} \circ \cdots \circ \Phi_{\xi(1) \xi(0)}
 \circ \pi_A(\xi)
\end{equation*}
 for any $\xi \in \Sigma_A$ and $m \geq 1$,
 we have
\begin{align*}
 d'(\pi_A(\xi),\pi_A(\xi'))
 \leq e^{-k\lambda'}
  \cdot d'(\pi_A \circ \sigma_A^k(\xi),
  \pi_A \circ \sigma_A^k(\xi'))
 \leq Ce^{\lambda'}d_\Sigma(\xi,\xi')^{\lambda'}
\end{align*}
  for any $\xi, \xi' \in \Sigma_A$,
 where $k$ is the integer with $d_\Sigma(\xi,\xi')=e^{-(k+1)}$.
\end{proof}

%
%
\subsection{Construction of a family of measures}
Put $\Delta= \{(p,q) \in M \times M \st q \in \cF^s_{\delta_0}(p)\}$.
We can define a continuous function $\eta$ on $\Delta$ 
 by $\eta(p,p)=0$
 and $\Phi^{\eta(p,q)}(p) \in \cF^{ss}(q)$.
Since the foliation $\cF^{ss}$ is H\"older continuous,
 so $\eta$ is.

Fix a positive H\"older continuous function $f$ on $M$.
We define a function $u$ on $\Delta$ by
\begin{displaymath}
u(p,q)=
 \int_0^\infty\left(
 f \circ \Phi^t(q)-f \circ \Phi^{t+\eta(p,q)}(p)
   \right) dt
 -\int_0^{\eta(p,q)}f \circ \Phi^t(p)\;dt.
\end{displaymath}
By the H\"older continuity of $f$
 and the inequality (\ref{eqn:Anosov F-ss}),
 the function $u$ is well-defined.
It is easy to check that the equations
\begin{align}
\label{eqn:u-1}
 u(p_1,p_3) &=u(p_1,p_2)+u(p_2,p_3)\\
\label{eqn:u-2}
 u(\Phi^t(p_1),\Phi^t(p_2))
 &= u(p_1,p_2) -\int_0^t f \circ \Phi^\tau(p_2) d\tau
 + \int_0^t f \circ \Phi^\tau(p_1) d\tau.
\end{align}
 hold for any $p_1,p_2,p_3 \in M$ and $t >0$
 when the both sides of the equations are well-defined.

\begin{lemma}
\label{lemma:u-rho}
The function $u$ is H\"older continuous.
\end{lemma}
\begin{proof}
Put
\begin{displaymath}
u_\infty(p,q)= 
 \int_0^\infty\left(
 f \circ \Phi^t(q)- f \circ \Phi^{t+\eta(p,q)}(p) 
 \right) dt.
\end{displaymath}
It is sufficient to show that $u_\infty$ is H\"older continuous.

Take positive constants $\delta$, $C$, $\beta_1$, and $\beta_2$
 such that
\begin{eqnarray}
\label{eqn:u-rho 1}
|f(p)-f(q)| &\leq & Cd(p,q)^\delta\\
\label{eqn:u-rho 2}
d(\Phi^t(p),\Phi^t(q)) &\leq & e^{\beta_1 t}d(p,q),\\
\label{eqn:u-rho 3}
d(\Phi^t(p),\Phi^t(p')) &\leq & Ce^{-\beta_2 t} 
\end{eqnarray}
 for any $p,q \in M$, $p' \in \cF_{\delta_0}^{ss}(p)$,
 and $t \geq 0$.
Fix $(p_1,q_1),(p_2,q_2) \in \Delta$.
Put $p'_i=\Phi^{\eta(p_i,q_i)}(p_i)$ and
 $T=-(2\delta \beta_1)^{-1}
  \log(d(p_1,p_2)^\delta+d(q_1,q_2)^\delta)$.
We may assume that $(p_1,q_1)$ and $(p_2,q_2)$
 are sufficiently close to each other so that
 $\Phi^T(q_i) \in \cF_{\delta_0}^{ss}(\Phi^T(p'_i))$ for $i=1,2$.
Then, we have
\begin{align*}
\lefteqn{|u_\infty(p_1,q_1)-u_\infty(p_2,q_2)|}\\
 & \leq 
 \int_0^T \left|f \circ \Phi^t(p'_1) -f \circ \Phi^t(p'_2)\right| dt
 +\int_0^T \left|f \circ \Phi^t(q_1) -f \circ \Phi^t(q_2)\right| dt\\
 & \quad + \int_T^\infty
  \left|f \circ \Phi^t(q_1) -f \circ \Phi^t(p_1')\right| dt
 +\int_T^\infty
  \left|f \circ \Phi^t(q_2) -f \circ \Phi^t(p_2')\right| dt\\
 & \leq 
 C(d(p'_1,p'_2)^\delta+d(q_1,q_2)^\delta)
 \int_0^T e^{\beta_1 \delta t} dt+
 2C^{1+\delta} \int_T^\infty e^{-\beta_2 \delta t} dt\\
 & \leq 
 \frac{C}{\beta_1\delta}
 (d(p'_1,p'_2)^\delta+d(q_1,q_2)^\delta)^{\frac{1}{2}}
 + \frac{C^{1+\delta}}{\beta_2 \delta}
 (d(p'_1,p'_2)^\delta+d(q_1,q_2)^\delta)^{\frac{\beta_2}{2\beta_1}}.
\end{align*}

Since $p'_i$ is H\"older continuous
 with respect to $(p_i,q_i)$,
 the function $u_\infty$ is H\"older continuous.
\end{proof}

We define a function $f_A$ on $\Sigma_A$ by
\begin{equation}
 f_A(\xi)=u(\pi_A(\sigma_A(\xi)),\pi_A(\xi)).
\end{equation}
Since $\sigma_A$, $\pi_A$, and $u$ are H\"older continuous,
 also $f_A$ is.
\begin{lemma}
\label{lemma:f-A} 
For any $p \in R_i$, any itinerary $(\tau_p,\xi_p)$
 of $p$, $m \geq 0$, and $t \in [\tau_p(m),\tau_p(m+1)]$,
\begin{equation}
\label{eqn:f-A}
\int_0^t f \circ \Phi^\tau(p) \;d\tau
 =u(\Phi^t(p),\pi_A(\sigma_A^m(\xi_p)))
 -u(p,\pi_A(\xi_p))
 +\sum_{j=0}^{m-1}f_A(\sigma_A^j(\xi_p)).
\end{equation}
\end{lemma}
\begin{proof}
Since $\Phi^{\tau}(p) \in \cF^s_{\delta_0}(\Phi^{\tau_p(j)}(p))$
 for any $\tau \in [\tau_p(j),\tau_p(j+1)]$, 
 we have $\eta(\Phi^t(p),\Phi^{t'}(p))=t-t'$,
 and hence,
\begin{equation*}
 u(\Phi^t(p),\Phi^{t'}(p))=\int_{t'}^t f \circ \Phi^\tau(p)\;d\tau
\end{equation*}
 for any $j=0,\cdots,m$ and $t,t' \in [\tau_p(j), \tau_p(j+1)]$.
Put $p_j=\pi_A \circ \sigma_A^j(\xi_p)$.
The equation (\ref{eqn:u-1}) implies 
\begin{align*}
 \int_{\tau_p(j)}^{\tau_p(j+1)} f \circ \Phi^\tau(p)\;d\tau
 & =u(\Phi^{\tau_p(j+1)}(p),p_{j+1})+u(p_{j+1},p_j)
  +u(p_j,\Phi^{\tau_p(j)}(p))\\
 & =u(\Phi^{\tau_p(j+1)}(p),p_{j+1})+f_A \circ\sigma_A^j(\xi_p)
  -u(\Phi^{\tau_p(j)}(p),p_j)
\end{align*}
 for any $j=0,\cdots,m-1$.
Remark that all terms in the above equation are well-defined 
 since $\{p_j, p_{j+1}\} \subset \cF^s_{\delta_0/2}(\Phi^\tau(p))$
 for any $\tau \in [\tau_p(j),\tau_p(j+1)]$.
Now, we have
\begin{align*}
\lefteqn{\int_0^t f \circ \Phi^\tau(p) \;d\tau}\\
 &= \int_{\tau_p(m)}^t f \circ \Phi^\tau(p) \;d\tau
 + \sum_{j=0}^{m-1}\int_{\tau_p(j)}^{\tau_p(j+1)}
  f \circ \Phi^\tau(p)\;d\tau\\
 &= u(\Phi^t(p),\Phi^{\tau_p(m)}(p))\\
 & \quad
  + \sum_{j=0}^{m-1}\left\{
  u(\Phi^{\tau_p(j+1)}(p),p_{j+1})+f_A(\sigma_A^j(\xi_p))
  -u(\Phi^{\tau_p(j)}(p),p_j)
 \right\}\\
 &= u(\Phi^t(p),\Phi^{\tau_p(m)}(p))+u(\Phi^{\tau_p(,m)}(p),p_m)
  -u(p,p_0)
 +\sum_{j=0}^{m-1}f_A(\sigma_A^j(\xi_p))\\
 & =u(\Phi^t(p),\pi_A(\sigma_A^m(\xi_p)))
 -u(p,\pi_A(\xi_p))
 +\sum_{j=0}^{m-1}f_A(\sigma_A^j(\xi_p)).
\end{align*}
\end{proof}

Since $u$ is bounded and $f$ is positive,
 the above lemma implies that
\begin{equation*}
\inf_{\xi \in \Sigma_A}\sum_{j=0}^{m-1}f_A \circ \sigma_A^{j}(\xi) >0
\end{equation*}
 for some $m \geq 1$.
By Lemma \ref{lemma:pressure},
 there exists $\rho_\Phi>0$ such that
 $P_{\sigma_A}(-\rho_\Phi \cdot f_A)=0$.
Let $(\mu, h)$ be the pair of a measure and a H\"older function
 on $\Sigma_A$ that is given by Theorem \ref{thm:Gibbs} for
 $-\rho_\Phi \cdot f_A$.
The measure $\mu$ satisfies
\begin{equation}
\label{eqn:mu}
\log \frac{d (\mu \circ \sigma_A)}{d\mu}(\xi)
 =\rho_\Phi \cdot f_A(\xi) 
\end{equation}
 for any $\xi \in \Sigma_A$.

Fix $\epsilon_1>0$ such that $B(p_1,4\epsilon_1) \subset \Int R_1$.
Put
\begin{align*}
\Sigma_A^* &=\left\{
 \xi \in \Sigma_A \st
 \cl{R}_{\xi(0)} \cap \pi_{\xi(0)}^{-1}(\pi_A(\xi))
 \subset {\bigcup}_{t>0}\Phi^{-t}(B(p_1,\epsilon_1))
\right\},\\
 \Sigma_A^\infty &=\bigcap_{m \geq 0}\sigma_A^{-m}(\Sigma_A^*).
\end{align*}
By the inclusion (\ref{eqn:Markov s}),
 we have $\sigma_A^{-1}(\Sigma_A^*) \subset \Sigma_A^*$,
 and hence, $\Sigma_A^\infty=\sigma_A(\Sigma_A^\infty)=
 \sigma_A^{-1}(\Sigma_A^\infty)$.
\begin{lemma}
\label{lemma:Sigma-*} 
The set $\Sigma_A^*$ has non-empty interior
 as a subset of $\Sigma_A$.
\end{lemma}
\begin{proof}
Since $\cF^s(p)= \bigcup_{t \in \RR}\Phi^t(\cF^{ss}(p))$ for any $p$,
 there exists $\delta'>0$ such that
\begin{equation}
\label{eqn:Sigma-* 1}
\cF^s_{\delta_0}(p) \subset
 \bigcup_{t \in (-\delta',\delta')} \Phi^t(\cF^{ss}_{\delta'}(p))
\end{equation}
 for any $p \in M$.
Since $\Phi$ contracts leaves of $\cF^{ss}$ uniformly,
 we can take $T'>0$ such that
\begin{equation}
\label{eqn:Sigma-* 2} 
\Phi^t(\cF^{ss}_{\delta'}(p)) \subset \cF^{ss}_{\epsilon_1}(\Phi^t(p))
\end{equation}
 for any $p \in M$ and $t \geq T'$.
Put $p_1'=\Phi^{-(T'+\delta')}(p_1)$
 and suppose $p_1' \in R_{i_1}$.
By the inclusions (\ref{eqn:Sigma-* 1}) and (\ref{eqn:Sigma-* 2}),
 we obtain
\begin{align*}
\cF^s_{\delta_0}(p_1')
&\subset
 \bigcup_{|t| <\delta'}\Phi^t(\cF^{ss}_{\delta'}(p_1'))\\
&\subset
 \bigcup_{|t| < \delta'}
 \Phi^{t-(T'+\delta')}(\cF^{ss}_{\epsilon_1}(p_1))
\subset
 \bigcup_{t>0}\Phi^{-t}(B(p_1,\epsilon_1)).
\end{align*}
Put $q_1=\pi_{i_1}(p_1')$.
Since $\cl{\pi_{i_1}^{-1}(q_1)} \subset \cF^s_{\delta_0}(p_1')$,
 there exists an open neighborhood $U$ of $q_1$
 in $\Lambda^u_{i_1}$ such that
\begin{equation*}
 \cl{\Int R_{i_1}} \cap \pi_{i_1}^{-1}(p)
 \subset \cl{\pi_{i_1}^{-1}(p)}
 \subset \bigcup_{t >0} \Phi^{-t}(B(p_1,\epsilon_1)).
\end{equation*}
 for any $p \in U$.
Therefore, $\pi_A^{-1}(U) \cap \{\xi \in \Sigma_A \st \xi(0)=i_1\}$
 is a non-empty open subset of $\Sigma_A^*$.
\end{proof}

\begin{lemma}
The restriction of $\pi_A$ to $\Sigma_A^\infty$ is injective and
\begin{equation}
 \mu(\Sigma_A^\infty)=1. 
\end{equation}
\end{lemma}
\begin{proof}
Fix $\xi \in \Sigma_A^\infty$
 and an itinerary $(\tau_p,\xi_p)$ of $p=\pi_A(\xi)$.
Since $\sigma_A^m(\xi) \in \Sigma_A^*$ for any $m \geq 1$,
 there exists a sequence $(t_k)_{k \geq 1}$ such that
 $\lim_{k \ra \infty}t_k= +\infty$ and
 $\Phi^{t_k}(p) \in B(p_1,\epsilon_1)$ for any $k \geq 1$.
By (\ref{eqn:Markov s}) and (\ref{eqn:Markov u}),
 the positive orbit $\{\Phi^t(p) \st t \geq 0\}$
 of $p$ intersects neither
 $\vphi_i([0,1] \times \del R_i^s \times R^u_i)$
 nor
 $\vphi_i([0,1] \times R_i^s \times \del R^u_i)$
 for any $i=1,\cdots,i_*$.
It implies that
 $\Phi^t(p) \in \Int R_{\xi_p(m)}$
 for any $m \geq 1$ and $t \in (\tau_p(m),\tau_p(m+1))$.
Since $\Int R_i \cap \Int R_j=\emptyset$
 for any $i \neq j$,
 $(\tau_p,\xi_p)$ is the unique itinerary of $p=\pi_A(\xi)$.
Therefore, the restriction of $\pi_A$ to $\Sigma_A^\infty$
 is injective.

Recall that $h$ is a positive H\"older continuous function
 and $h \cdot \mu$ is a $\sigma_A$-invariant ergodic measure.
Since $\sigma_A^{-1}(\Sigma_A^*) \subset \Sigma_A^*$,
\begin{equation}
 h \cdot \mu(\Sigma_A^\infty)
  =\lim_{m \ra \infty} h \cdot \mu(\sigma_A^{-m}(\Sigma_A^*))
 =\mu(\Sigma_A^*).
\end{equation} 
Lemma \ref{lemma:Sigma-*} implies $\mu(\Sigma_A^*)>0$.
Since $h \cdot \mu$ is ergodic, we have $\mu(\Sigma_A^\infty)=1$.
\end{proof}

We define a probability measure $\mu'$ on $\Lambda^u$
 by the pullback of $\mu$ by $\pi_A$,
 that is, $\mu'(B)=\mu(\pi_A^{-1}(B))$
 for any Borel subset $B$ of $\Lambda^u$.
Recall that $\mu$ is non-atomic and
 positive on each open subset of $\Sigma_A$.
Since $\pi_A$ is a surjective map onto $\Lambda^u$,
 the support of $\mu'$ is $\Lambda^u$.
By the above lemma,
 we have $\mu'(\Lambda^u \-- \pi_A(\Sigma_A^\infty))=0$.
Since the restriction of $\pi_A$ on $\Sigma_A^\infty$ is injective,
 $\mu'$ is non-atomic.

Take $\epsilon_*>0$ such that
\begin{equation}
\label{eqn:e-*}
 \Phi^t(\cF^s_{\epsilon_*}(p)) \subset \cF^s_{\epsilon_1}(\Phi^t(p)) 
\end{equation}
 for any $p \in M$ and $t \geq 0$.
We put
\begin{equation*}
V_i = \bigcup_{p \in \cl{R}_i} \cF^s_{\epsilon_*}(p),\hsp
V_i^\infty
 =V_i \cap \pi_i^{-1}(\pi_A(\Sigma_A^\infty)).
\end{equation*}
Remark that $V_i^\infty$ is a subset of
 $\bigcup_{t>0} \Phi^{-t}(B(p_1,2\epsilon_1))$.

Recall that $\pi_i$ is the projection from a neighborhood $U_i$
 of $\cl{R_i}$ to a neighborhood $I^u_i$ of $\Lambda^u_i$.
For $p \in M$,
 let $h_{i,p}$ be the restriction of $\pi_i$ to $\cF^{uu}(p) \cap U_i$.
Since $h_{i,p}$ is a covering map,
 we can define a measure $\nu_{i,p}$ on $\cF^{uu}(p)$ by
\begin{equation}
\label{eqn:def of nu}
 \nu_{i,p}(A)
 =\rho_\Phi \cdot \int_{A \cap V_i}
  \exp(u(q,\pi_i(q))) \;d(\mu' \circ h_{i,p})(q)
\end{equation}
  for any Borel subset $A$ of $\cF^{uu}_i(p)$.
The measure $\nu_{i,p}$ is non-atomic, locally finite,
 and positive
 on each non-empty open subset of $\cF^{uu}(p)\cap V_i$.
Moreover, it satisfies $\nu_{i,p}(\cF^{uu}(p)\-- V_i^\infty)=0$.
By the equation (\ref{eqn:mu}), we have
 \begin{equation}
 \label{eqn:mu nu}
 \log\frac{d(\nu_{\xi(1),\pi_A(\sigma_A(\xi))}
   \circ \Phi_{\xi(1)\xi(0)})}{d\nu_{\xi(0),\pi_A(\xi)}}(\xi_A(\xi)) 
 = \rho_\Phi \cdot f_A(\xi)
\end{equation}
 for any $\xi \in \Sigma_A$.

The following proposition is a keystone to 
 the construction of the required family $\{\nu_p\}_{p \in M}$.
\begin{prop}
\label{lemma:nu 1}
The equation
\begin{equation}
\label{eqn:nu 1}
\log
 \frac{d(\nu_{\xi_p(m),\Phi^t(q)} \circ \Phi^t)}{d\nu_{\xi_p(0),q}}(q)
 =\rho_\Phi \cdot \int_0^t f \circ \Phi^\tau(q)\;d\tau.
\end{equation}
 holds for any $p \in M$, any itinerary $(\tau_p,\xi_p)$ of $p$,
 $m \geq 1$, $t \in [\tau_p(m),\tau_p(m+1)]$,
 and $q \in \cF^s_{\epsilon_*}(p)$.
\end{prop}
\begin{proof}
Put $q'=\Phi^t(q)$, $\xi_p(j)=i_j$,
 and $p_j=\pi_A(\sigma_A^j(\xi_p))$ for $j=0,\cdots,m$.
Notice that $q' \in
 \cF^s_{\epsilon_1}(\Phi^t(p)) \subset \cF^s_{\delta_0}(p_m)$.
On a small neighborhood of $q$ in $\cF^{uu}(q)$, we have
\begin{equation*}
\Phi^t=h_{q'p_m}
 \circ \Phi_{i_m i_{m-1}} \circ \cdots
 \circ \Phi_{i_1 i_0} \circ h_{p_0 q}.
\end{equation*}
By the definition (\ref{eqn:def of nu}),
 we have
\begin{equation*}
\frac{d \nu_{i,q} \circ h_{qp}}{d \nu_{i,p}}
=\exp(\rho_\Phi \cdot u(q,\pi_i(q))-\rho_\Phi \cdot (p,\pi_i(p)))
=\exp(\rho_\Phi \cdot u(q,p)).
\end{equation*}
 if $\pi_i(p)=\pi_i(q)$.
It implies that
\begin{align*}
\lefteqn{
 \frac{d(\nu_{i_m,q'} \circ \Phi^t)}{d\nu_{i_0,q}}(q)}\\
&=
 \frac{d(\nu_{i_m,q'} \circ h_{q'p_m})}{d \nu_{i_m,p_m}}(p_m)
 \cdot
 \prod_{j=0}^{m-1}
  \frac{d(\nu_{i_{j+1},p_{j+1}} \circ \Phi_{i_{j+1}i_j})}{
   d\nu_{i_j,p_j}}(p_j)
 \cdot
 \frac{d(\nu_{i_0,p_0} \circ h_{p_0 p})}{d\nu_{i_0,q}}(q)\\
&=\exp\left(\rho_\Phi
 \cdot \left\{u(q',p_m)+\sum_{j=0}^{m-1}f_A(\sigma_A^j(\xi_p))
 +u(p_0,q)\right\}\right).
\end{align*}
By Lemma \ref{lemma:f-A} and the equation (\ref{eqn:u-1}),
 we obtain the required equation (\ref{eqn:nu 1}).
\end{proof}

\begin{cor}
\label{cor:nu}
For any $i=1,\cdots,i_*$, $q \in V_i^\infty$,
 and $t>0$ with $\Phi^t(q) \in B(p_1,2\epsilon_1)$,
\begin{equation*}
\frac{d(\nu_{1,\Phi^t(q)} \circ \Phi^t)}{d \nu_{i,q}}
 =\rho_\Phi \cdot \int_0^t f \circ \Phi^t(q) \;dt
\end{equation*}
In particular,
 $\nu_{i,p}=\nu_{j,p}$ on $\cF^{uu}(p) \cap V_i \cap V_j$
 for any $p \in M$.
\end{cor}
\begin{proof}
For $q \in V_i^\infty$,
 there exists $p \in \cF^s_{\epsilon_*}(q) \cap \cl{R}_i \cap V_i^\infty$.
Let $(\tau_p,\xi_p)$ be a itinerary of $p$ with $\xi_p(0)=i$
 and $m$ be the integer such that $t \in [\tau_p(m), \tau_p(m+1))$.
Since $\Phi^t(p) \in \cF^s_{\epsilon_1}(\Phi^t(q))
 \subset B(p_1,3\epsilon_1)$,
 we have $\tau_p(m)=1$.
Now, the former half of the corollary follows
 from the proposition immediately.

For any $p \in V_i^\infty$,
 there exists $t >0$ such that $\Phi^t(p) \in B(p_1,2\epsilon_1)$.
Hence, the former half implies $\nu_{i,p}=\nu_{j,p}$ on
 $\cF^{uu}(p) \cap V_i^\infty \cap V_j^\infty$
 for any $i,j=1,\cdots,i_*$ and $p \in M$.
Since $\nu_{i,p}(\cF^{uu}(p) \--V_i^\infty)=0$,
 we obtain the latter half of the corollary.
\end{proof}

Now, we prove Theorem \ref{thm:Radon}.
Recall that $M=\bigcup_{i=1}^{i_*} V_i$.
By Corollary \ref{cor:nu},
 we can define a measure $\nu_p$ on $\cF^{uu}(p)$
 by $\nu_p(U \cap V_i)=\nu_{i,p}(U \cap V_i)$
 for any $i=1,\cdots,i_*$ and any Borel subset $U$ of $\cF^{uu}(p)$.
Then, $\nu_p$ is a non-atomic locally finite measure
 whose support is $\cF^{uu}(p)$.

Fix $p \in M$
 and put $S=\{i=1,\cdots,i_* \st p \in \cl{R}_i\}$.
Notice that $\bigcup_{i \in S} \cl{R}_i$ is a neighborhood of $p$.
Take an open interval $I$ in
 $\cF^{uu}(p) \cap \bigcup_{i \in S}\cl{R}_i$ that contains $p$.
If $q \in M$ is sufficiently close to $p$,
 then $h_{qp}(I) \subset \bigcup_{i \in S} V_i$
 and $h_{qp}(I \cap V_i^\infty)=h_{qp}(I) \cap V_i^\infty$
 for any $i \in S$.
For any $i \in S$
 and any Borel subset $U$ of $I$, we have
\begin{align*}
\nu_{i,q}(h_{qp}(U))
&= \nu_{i,q}(h_{qp}(U) \cap V_i^\infty)\\
&= \nu_{i,q}(h_{qp}(U \cap V_i^\infty))\\
&=\int_{U \cap V_i^\infty} \exp\left(
 \rho_\Phi \cdot u(h_{qp}(p'),p')\right)d \nu_{i,p}(p')\\
&=\int_{U} \exp\left(
 \rho_\Phi \cdot u(h_{qp}(p'),p')\right)d \nu_{i,p}(p'),
\end{align*}
 and hence,
\begin{equation}
\label{eqn:nu holonomy}
\log \frac{d(\nu_q \circ h_{qp})}{d\nu_p}(p)
 =\rho_\Phi \cdot u(h_{qp}(p),p).
\end{equation}
Since $u$ is H\"older continuous,
 the family $\{\nu_p\}_{p \in M}$ is an element of $\cM(\Phi)$.
Since the manifold $M$ is closed and
\begin{equation*}
u(\Phi^t(p),p)=\int_0^t f \circ \Phi^\tau(p)\;d\tau 
\end{equation*}
 for any $p \in M$ and any sufficiently small $t>0$,
 the equation (\ref{eqn:nu holonomy}) implies
\begin{equation*}
\log \frac{d(\nu_{\Phi^t(p)} \circ \Phi^t)}{d\nu_p}(p)
 =\rho_\Phi \cdot \int_0^t f\circ \Phi^\tau(p) \;d\tau
\end{equation*}
 for any $p \in M$ and any $t \in \RR$.

\section{Deformation of Anosov flows}
\label{section:main thm}
\subsection{$C^{r+,0+}$-flows tangent to a foliation}
\label{sec:tangent}
Let $V$ be an open subset of
 $\RR^{n+1}=\{(t,x,y) \in \RR \times \RR^{n-1} \times \RR\}$.
A continuous map $F:V \ra \RR^n$ is called {\it of class $C^{1+,0+}$}
 if $F$, $\del F/\del x_1,\cdots, \del F/\del x_{n-1}$,
 and $\del F/ \del t$ are well-defined and H\"older continuous.
We say a sequence $\{F_k\}_{k \geq 1}$ of $C^{1+,0+}$ maps
 from $V$ to $\RR^n$ {\it converges} to $F$
 {\it in the $C^{1,0}$-topology}
 as $k \ra  \infty$
 if $F_k$,
 $\del F_k/\del x_1,\cdots, \del F_k/\del x_{n-1}$,
 and $\del F_k/ \del t$ converges to
 $F$, $\del F/\del x_1,\cdots, \del F/\del x_{n-1}$,
 and $\del F/ \del t$ with respect to the $C^0$-norm, respectively.

Let $M$ be a $C^{1+}$ $n$-dimensional closed manifold
 and $\cF$ a $C^{1+}$ codimension-one foliation on $M$.
A continuous flow $\Phi=\{\Phi^t\}_{t \in \RR}$
 is called {\it of class $C^{1+,0+}_{\cF}$} if
 it preserves each leaf of $\cF$ and
 there exists a family $\{\vphi_p:U_p \ra \RR^n\}_{p \in M}$
 of $C^{1+}$ foliation coordinates of $\cF$
 such that $U_p$ contains $p$ and the map
\begin{equation*}
 F_p:(w,t) \mapsto \vphi_p \circ \Phi^t \circ \vphi_p^{-1}(w)
\end{equation*}
 is well-defined and of class $C^{1+,0+}$
 on a neighborhood $V_p$ of $(\vphi_p(p),0)$ in $\RR^{n+1}$
 for any $p \in M$.
If $\Phi$ is of class $C^{1+,0+}_{\cF}$
 then for any $C^{1+}$ foliation coordinate $\vphi_p$
 and any sufficiently small neighborhood $V_p$ of $(\vphi_p(p),0)$,
 the above map $F_p$ is well-defined
 and of class $C^{1+,0+}$ on $V_p$.
Remark that any $C^{1+}$ flow which preserves each leaf of $\cF$
 is a $C^{1+,0+}_{\cF}$ flow.

We say a sequence $\{\Phi_k\}_{k \geq 1}$
 of $C^{1+,0+}_{\cF}$ flows {\it converges} to $\Phi$
 {\it in the $C^{1,0}$}-topology
 if there exists a neighborhood $V'_p \subset V_p$ of
 $(\vphi_p(p),0)$ such that
 the map
 $F_{k,p}:
  (w,t) \mapsto \vphi_p \circ \Phi_k^t \circ \vphi_p^{-1}(w)$
 is well-defined on $V'_p$ for any $k \geq 1$
 and the sequence $(F_{k,p})_{k \geq 1}$
 converges to $F$ in the $C^{1,0}$-topology.

Fix a H\"older continuous metric $g$ on $M$.
Let $\Phi$ be a $C^{1+,0+}_{\cF}$ flow on $M$.
It is generated by a H\"older continuous vector field tangent to $T\cF$.
Let $TM=T\Phi \oplus E_\Phi \oplus E_\Phi^\perp$
 be the orthogonal splitting with respect to $g$
 that satisfies $T\Phi \oplus E_\Phi= T\cF$.
By $\pi_\Phi$ and $\pi_\Phi^\perp$,
 we denote the natural projection onto $E_\Phi$ and $E_\Phi^\perp$
 with respect to the splitting.
For each $t \in \RR$,
 the restriction of $\Phi^t$ to each leaf of $\cF$
 defines an isomorphism $D_{\cF}\Phi^t$ of $T\cF$.
The one-parameter family $\{D_\cF\Phi^t\}_{t \in \RR}$
 forms a H\"older continuous flow.
We define a H\"older continuous flow
 $N_{\cF}\Phi=\{N_\cF\Phi^t\}_{t \in\RR}$ on $E_\Phi$ by
\begin{equation*}
  N_{\cF}\Phi^t=\pi_\Phi \circ D_{\cF} \Phi^t|_{E_\Phi}.
\end{equation*}
Linear holonomy of $\cF$ along the orbits of $\Phi$
 also induces a H\"older continuous flow $N_\cF^\perp \Phi$
 on $E_\Phi^\perp$.

For a flow $\psi=\{\psi^t\}_{t \in \RR}$ on a topological space $B$,
 we say a function $\alpha$ on $B \times \RR$ is
 {\it a cocycle} over a flow $\psi$
 if $\alpha(p,s+t)=\alpha(p,s)+\alpha(\psi^s(p),t)$
 for any $p \in B$ and $s,t \in \RR$.
Put
\begin{align*}
\alpha_{\Phi}(p,t;\cF;g) &
 =\log {\det}_g(N_\cF\Phi^t)_p,\\
\alpha_{\Phi}^\perp(p,t;\cF;g) &
 =\log {\det}_g(N_\cF^\perp \Phi^t)_p,\\
\beta_{\Phi}(p,t;\cF;g) &
 =\log \|(N_\cF\Phi^t)_p\|_g.
\end{align*}
Then, $\alpha_\Phi$ and $\alpha_\Phi^\perp$ are cocycles over $\Phi$,
 and $\beta_\Phi$ satisfies
\begin{equation*}
 \beta(p,s+t) \leq \beta(p,s)+\beta(\Phi^s(p),t).
\end{equation*}
If a sequence $(\Phi_k)_{k \geq 1}$ of $C^{1+,0+}_\cF$
 flows on $M$ converges to $\Phi$ in the $C^{1,0}$-topology,
 then $\alpha_{\Phi_k}$, $\alpha_{\Phi_k}^\perp$,
 and $\beta_{\Phi_k}$ converges uniformly to
 $\alpha_{\Phi}$, $\alpha_{\Phi}^\perp$,
 and $\beta_{\Phi}$ respectively on $M \times [-T,T]$
 for any given $T>0$.

\begin{lemma}
\label{lemma:Anosov}
Let $\Phi$ be a $C^{1+}$ and $C^{1+,0+}_{\cF}$ flow on $M$.
Suppose
\begin{equation*}
\inf_{p \in M}
 \left\{\alpha_{\Phi}^\perp(p,T_1;\cF,g),
  -\beta_{\Phi}(p,T_1;\cF,g)\right\}>0
\end{equation*}
  for some $T_1>0$.
Then, there exists an Anosov splitting
 $TM=T\Phi \oplus E^{ss} \oplus E^{uu}$
 such that $T\Phi \oplus E^{ss}=T\cF$.

Moreover, if $\Phi$ is of class $C^\infty$ and
\begin{equation*}
\inf_{p \in M} \left\{
 \alpha_{\Phi}(p,T_2;\cF,g)
   +\alpha_{\Phi}^\perp(p,T_2;\cF,g)
   -\beta_{\Phi}(p,T_2;\cF,g)
 \right\}>0
\end{equation*}
  for some $T_2>0$ in addition,
 then the subbundle $T\Phi \oplus E^{uu}$ is of class $C^{1+}$.
\end{lemma}
\begin{proof}
Without loss of generality,
 we may assume that $T\cF$ is transversely orientable.
Let $E$ be the orthogonal complement of $T\Phi$ with respect to $g$
 and $\pi^\perp$ the orthogonal projection from $TM$ to $E$.
We define a flow $N\Phi=\{N\Phi^t\}_{t \in \RR}$ on $E$ by
 $N\Phi^t=\pi^\perp \circ D\Phi^t$.
Remark that $E=E_\Phi \oplus E_\Phi^\perp$
 and the projection $\pi^\perp$ coincides with $\pi_\Phi$
 on $T\Phi \oplus E_\Phi$.
In particular, $N\Phi^t(v)=N_\cF\Phi^t(v)$ for any $v \in E_\Phi$.
We also remark that $N_\cF^\perp\Phi^t(v')
 = \pi_\Phi^\perp \circ D\Phi^t(v')
 = \pi_\Phi^\perp \circ N\Phi^t(v')$ for any $v' \in E_\Phi^\perp$.

Let $\Gamma(E_\Phi)$ be the set of continuous sections of
 $E_\Phi$.
It becomes a Banach space by a norm
 $\|\xi\|_\Gamma=\sup_{p \in M}\|\xi(p)\|_g$.
Let $v^\perp$ be the unit tangent vector field of $E_\Phi^\perp$.
Since
\begin{equation*}
\pi_\Phi^\perp \circ N\Phi^t(v^\perp(p))
=N_\cF^\perp\Phi^t(v^\perp(p))
=\exp(\alpha_\Phi^\perp(p,t;\cF,g))  \cdot v^\perp(\Phi^t(p)),
\end{equation*}
 we can define a map $\Gamma\Phi^t$ on $\Gamma(E_\Phi)$ by
\begin{equation*}
\Gamma\Phi^t(\xi)(\Phi^t(p))+v^\perp(\Phi^t(p))
 = \exp(-\alpha_\Phi^\perp(p,t;\cF,g)) \cdot 
  N\Phi^t\left(\xi(p)+v^\perp(p)\right). 
\end{equation*}
Since $N\Phi=\{N\Phi^t\}$ is a flow
 and $\alpha_\Phi$ is a cocycle,
 the family $\{\Gamma\Phi^t\}_{t \in \RR}$
 forms a flow on $\Gamma(E_\Phi)$.
Put
\begin{equation*}
 C_\Phi(t)=\inf_{p \in M}\left(
 \alpha_\Phi^\perp(p,t;\cF,g)-\beta_\Phi(p,t;\cF,g)\right).
\end{equation*}
Then, we have
 $C_\Phi(t+t') \geq C_\Phi(t)+C_\Phi(t')$ and
\begin{equation}
\label{eqn:Anosov-1}
\|\Gamma\Phi^t(\xi)-\Gamma\Phi^t(\xi')\|_\Gamma
  \leq \exp(-C_\Phi(t))\|\xi-\xi'\|_\Gamma
\end{equation}
 for any $t,t' \geq 0$
 and $\xi,\xi' \in \Gamma(E_\Phi^\perp)$.
Since $C_\Phi(T_1)>0$,
 the contraction mapping principle implies that
 the flow $\Gamma\Phi$ admits a unique fixed point
 $\xi_0$.
Let $E_0$ be the one-dimensional subbundle
 of $TM$ that is generated by the vector field $v^\perp+\xi_0$.
Since the restriction of $\pi_\Phi^\perp$ to $E_0$
 is an isomorphism onto $E_\Phi^\perp$,
 there exists $T'>0$ such that
 $\inf_{p \in M}\|N\Phi^{T'}|_{E_0(p)}\|_g >1$.
By the standard argument
 (see {\it e.g.} \cite{Do}),
 the flow $\Phi$ admits an Anosov splitting
 $TM=T\Phi \oplus E^{ss} \oplus E^{uu}$
 with $T\Phi \oplus E_\Phi=T\Phi \oplus E^{ss}$
 and $T\Phi \oplus E_0=T\Phi \oplus E^{uu}$.

The latter half of the lemma
 is a consequence of the $C^r$-section theorem for $E_0$
 (see {\it e.g.} \cite{Sh}).
\end{proof}

\subsection{Replacement of smooth structures}
Let $M$ be a $C^\infty$ closed manifold of dimension $n$
 and $\Phi$ be a $C^\infty$ codimension-one Anosov flow on $M$.
Let $TM=T\Phi \oplus E^{ss} \oplus E^{uu}$ be
 the Anosov splitting,
 $\cF$ be the weak stable foliation,
 and $\cF^{uu}$ be the strong unstable foliation of $\Phi$.
Recall that $\cF$ is a $C^{1+}$ foliation
 with $C^\infty$ leaves.

Fix an integer $r \geq 1$
 and a $C^{r+}$ one-dimensional foliation $\cL$
 that is transverse to $\cF$.
By $\pi_x$ and $\pi_y$,
 we denote the natural projections from
 $\RR^n=\RR^{n-1} \times \RR$ to $\RR^{n-1}$ and $\RR$, respectively.
For each $p \in M$,
 we can take a $C^{1+}$ coordinate $\vphi_p:U_p \ra \RR^n$
 and a $C^{1+}$ map $L_p:(-1,1) \ra \cF^{uu}(p)$
 that satisfy the following conditions:
\begin{enumerate}
 \item $\vphi_p(p)=(0,0)$ and $\vphi_p(U_p)=(-1,1)^n$.
 \item $\vphi_p^{-1}((-1,1)^{n-1} \times y) \subset \cF(q)$
 and $\vphi_p^{-1}(x \times (-1,1)) \subset \cL(q)$
 for any $q =\vphi_p^{-1}(x,y) \in U_p$.
 \item $\pi_x \circ \vphi_p$ is of class $C^{r+}$.
 \item $L_p(0)=p$ and
 $L_p(y) \in \cF^{uu}(p) \cap \vphi_p^{-1}((-1,1)^{n-1} \times y)$
 for any $y \in (-1,1)$.
\end{enumerate}
For $p,q \in M$, 
 there exist $V_{pq}^x \subset (-1,1)^{n-1}$,
 $V_{pq}^y \subset (-1,1)$,
 and maps $\vphi_{pq}^x:V_{pq}^x \ra (-1,1)^{n-1}$
 and $\vphi_{pq}^y:V_{pq}^y \ra (-1,1)$ such that
\begin{align*}
\vphi_q(U_p \cap U_q)
 &= V_{pq}^x \times V_{pq}^y,\\
\vphi_p \circ \vphi_q^{-1}(x,y)
 &= (\vphi_{pq}^x(x),\vphi_{pq}^y(y)).
\end{align*}

Since $\Phi$ is Anosov,
 we can take a $C^\infty$ metric $g$ on $TM$
 and $\lambda>0$ such that
\begin{equation}
\label{eqn:Phi Anosov}
\inf_{p \in M}\left\{
 \alpha_\Phi^\perp(p,t;\cF,g),
 -\beta_\Phi(p,t;\cF,g)
\right\} \geq \lambda \cdot t
\end{equation}
 for any $t >0$.
Put $f(p)=-(\del \alpha_\Phi/\del t)(p,0;\cF,g)$.
Since $T\cF$ is of class $C^{1+}$
 and $\alpha_\Phi \leq \beta_\Phi$,
 the function $f$ is H\"older continuous and positive.
By Theorem \ref{thm:Radon},
 there exists an element $\nu=\{\nu_p\}_{p \in M}$ of $\cM(\Phi)$
 and a constant $\rho_\Phi>0$ such that
\begin{equation}
\label{eqn:nu rho}
\frac{d(\nu_{\Phi^t(p)} \circ \Phi^t)}{d\nu_p}(p)
 =\rho_\Phi \cdot \int_0^t f \circ \Phi^\tau(p) d\tau
 = -\rho_\Phi \cdot \alpha_\Phi(p,t;\cF,g)
\end{equation}
 for any $p \in M$ and $t \in \RR$.
Since $L_p$ is a homeomorphism onto its image,
 the measure $\nu_p \circ L_p$ on $(-1,1)$ is well-defined
 for each $p \in M$.
We define a function $\eta_p$ on $(-1,1)$ by
\begin{equation*}
\eta_p(y)=\int_0^y d(\nu_p \circ L_p).
\end{equation*}

\begin{lemma}
\label{lemma:Holder}
For each $p \in M$,
  $\eta_p$ is a bi-H\"older homeomorphism onto its image.
\end{lemma}
\begin{proof}
Since $\nu_p$ is non-atomic, locally finite, and positive
 on each non-empty open subset of $\cF^{uu}(p)$,
 the map $\eta_p$ is a homeomorphism onto its image.

The metric $g$ on $M$ induces
 the Lebesgue measure $\nu^L_p$ on $\cF^{uu}(p)$ for each $p \in M$.
There exist constants $C>0$ and $\lambda_+>\lambda_->0$ such that
 $e^{\lambda_- t -C} \leq \|D\Phi^t|_{E^{uu}(p)}\| \leq
 e^{\lambda_+ t + C}$
 for any $p \in M$ and any $t > 0$.
Take $\lambda_+'>\lambda_-'>0$ such that
 $\lambda_-' \leq \rho_\Phi \cdot f(p) \leq \lambda_+'$
 for any $p \in M$.
Then, we have
\begin{align*}
 e^{\lambda_- t-C} \cdot \nu^L_p(U) &\leq \nu^L_{\Phi^t(p)}(U)
 \leq e^{\lambda_+ t+C} \cdot \nu^L_p(U)\\
 e^{\lambda_-' t} \cdot \nu_p(U) &\leq \nu_{\Phi^t(p)}(U)
 \leq e^{\lambda_+' t} \cdot \nu_p(U)
\end{align*}
 for any $p \in M$, any $t \geq 0$,
 and any Borel subset $U$ of $\cF^{uu}(p)$.
Since the function $(p,q) \mapsto (d(\nu_p \circ h_{pq})/d\nu_q)(q)$
 is continuous,
 there exists $C'>0$ such that
 $e^{-C'} \leq \nu_p(I_p) \leq e^{-C'}$
 for any $p \in M$ and any interval $I_p$ in $\cF^{uu}(p)$
 with $\nu^L_p(I_p)=1$.
For any $p \in M$ and any interval $I$ in $\cF^{uu}(p)$
 with $\nu^L_p(I) < 1$,
 there exists $t \geq 0$ such that $\nu^L_p(\Phi^t(I))=1$.
We have
\begin{align*}
e^{-\lambda_+ t -C}
  & \leq \nu^L_p(I) \leq e^{-\lambda_- t +C},\\
e^{-\lambda_+' t -C'} & \leq \nu_p(I) \leq e^{-\lambda_-' t +C'},
\end{align*}
 and hence, 
\begin{align*}
\nu_p(I)^{\lambda_+} &
 \leq e^{\lambda_+ C' + \lambda_-' C} \cdot \nu^L_p(I)^{\lambda_-'},\\
\nu^L_p(I)^{\lambda_+'} &
 \leq e^{\lambda_- C' + \lambda_+' C} \cdot \nu_p(I)^{\lambda_-}.
\end{align*}
Since $(\vphi_p^{-1})^*g$ is comparable with
 the Euclidean metric on $\vphi_p(U_p)$,
 the above inequalities shows that $\nu_p$ is bi-H\"older continuous.
\end{proof}

\begin{lemma}
\label{lemma:ch-vphi}
For any $p,q \in M$,
 the map $\eta_p \circ \vphi_{pq}^y \circ \eta_q^{-1}$
 on $\eta_p(V_{pq}^y)$ is of class $C^{1+}$.
\end{lemma}
\begin{proof}
Fix $p \in M$ and $y_0 \in (-1,1)$.
Suppose that $q \in M$ is sufficiently close to $p$
 so that the holonomy map
 $h_{pq}:\cF^{uu}(q) \ra \cF^{uu}(p)$
 of $\cF^s$ is well-defined on a neighborhood of $L_p(y_0)$.
Since $\vphi_{pq}^y=L_p^{-1} \circ h_{pq} \circ L_q$,
 we have
\begin{equation*}
\frac{d (\nu_p \circ L_p \circ \vphi_{pq}^y)}{d (\nu_q \circ L_q)} 
 =\frac{d (\nu_p \circ h_{pq})}{d \nu_q} \circ L_q.
\end{equation*}
Therefore,
\begin{align*}
\lefteqn{\eta_p \circ \vphi_{pq}^y \circ \eta_q^{-1}(y)
 -\eta_p \circ \vphi_{pq}^y \circ \eta_q^{-1}(y_0)}\\
 & =\int_{\vphi_{pq}^y \circ \eta_q^{-1}(y_0)
 }^{\vphi_{pq}^y \circ \eta_q^{-1}(y)} d(\nu_p \circ L_p)\\
 & = \int_{\eta_q^{-1}(y_0)}^{\eta_q^{-1}(y)}
  \left(\frac{d(\nu_p \circ h_{pq})}{d \nu_q} \circ L_q\right)
  \; d(\nu_q \circ L_q)\\
 &= \int_{y_0}^y 
  \frac{d(\nu_p \circ h_{pq})}{d \nu_q} \circ
 (L_q \circ \eta_q^{-1})(\eta)
  d\eta.
\end{align*}
Since $d(\nu_p \circ h_{pq})/d \nu_q$ is H\"older continuous
 by the definition of $\cM(\Phi)$,
 the map $\eta_p \circ \vphi_{pq}^y \circ \eta_q^{-1}$
 is of class $C^{1+}$ at $y_0$.
\end{proof}

We define a map $\ch{\vphi}_p:U_p \ra \RR^n$ by
\begin{equation*}
 \ch{\vphi}_p \circ \vphi_p^{-1}(x,y)=(x,\eta_p(y)).
\end{equation*}
The above lemma implies that
 the family $\{\ch{\vphi}_p\}_{p \in M}$ defines a $C^{1+}$-structure
 on $M$.
We denote this $C^{1+}$-manifold by $\ch{M}$.
The identity map on $M$ as a set
 induces a map $i_M:M \ra \ch{M}$.
It is bi-H\"older by Lemma \ref{lemma:Holder},
 but is not of class $C^{1+}$ in general.
Put $\ch{\cF}=i_M(\cF)$ and $\ch{\cL}=i_M(\cL)$.
By the definitions of $\vphi_p$ and $\ch{\vphi}_p$,
 both $\ch{\cF}$ and $\ch{\cL}$ are $C^{1+}$ foliations
 and the restriction of $i_M$ to each leaf of $\cF$ is
 of class $C^{1+}$.
Hence, the map $i_M$ induces a bi-H\"older isomorphism
 $D_x i_M$ from $T\cF$ to $T\ch{\cF}$.

We define a continuous flow $\ch{\Phi}=\{\ch{\Phi}^t\}_{t \in \RR}$
 on $\ch{M}$ by $\ch{\Phi}^t= i_M \circ \Phi^t \circ i_M^{-1}$.
Since
\begin{equation*}
\ch{\vphi}_p \circ \ch{\Phi}^t \circ \ch{\vphi}_p^{-1}
 =(\Id \times \eta_p)
 \circ (\vphi_p \circ \Phi^t \circ \vphi_p^{-1})
 \circ (\Id \times \eta_p^{-1}),
\end{equation*}
 it is a $C^{1+,0+}_{\ch{\cF}}$ flow
 by Lemma \ref{lemma:ch-vphi}.
The flow
 $\ch{\Phi}$ is generated by a H\"older continuous vector field.
We denote it by $\ch{X}$.
Take a H\"older continuous metric $\ch{g}$ such that
 $T_p \ch{M}=T\ch{\cF}(p) \oplus \RR \cdot D\ch{\vphi}_p(\del/\del y)$
 is an orthogonal splitting with respect to $\ch{g}$
 for any $p \in \ch{M}$,
 $(D_x i_M)^*(\ch{g}|_{T\ch{\cF}}) = g|_{T\cF}$,
 and $\|D\ch{\vphi}_p^{-1}(\del/\del y)\|_{\ch{g}}=1$.
We can check that
\begin{align*}
\alpha_{\ch{\Phi}}(i_M(p),t;\ch{\cF},\ch{g}) &
 =\alpha_{\Phi}(p,t;\cF,g),\\
\alpha_{\ch{\Phi}}^\perp(i_M(p),t;\ch{\cF},\ch{g}) &
 = \frac{d(\nu_{\Phi^t(p)} \circ \Phi^t)}{d \nu_p}(p)
 = -\rho_{\Phi} \cdot \alpha_{\Phi}(p,t;\cF,g),\\
\beta_{\ch{\Phi}}(i_M(p),t;\ch{\cF},\ch{g}) &
 =\beta_{\Phi}(p,t;\cF,g).
\end{align*}

\begin{lemma}
\label{lemma:stability} 
Let $(\ch{\Phi}_k)_{k \geq 1}$ be a sequence of $C^{1+,0+}$ flows
 on $\ch{M}$ that converges to $\ch{\Phi}$ in the $C^{1,0}$ topology.
Then, $\ch{\Phi}_k$ is topologically equivalent to $\Phi$
 for any sufficiently large $k$.
\end{lemma}
\begin{proof}
We mimic the proof of the structural stability of Anosov flows.
Let $\cF^{ss}$ be the strong stable foliation of $\Phi$.
Put $\ch{\cF}^{ss}=i_M(\cF^{ss})$.
It is a H\"older foliation and
 its restriction to each leaf of $\ch{\cF}$ is of class $C^{1+}$.
Let $d^{ss}_p$ be the distance on $\ch{\cF}^{ss}(p)$
 that is induced from the restriction of the metric
 $\ch{g}$ to $\ch{\cF}^{ss}(p)$.
Since $D_x i_M$ is an isometry,
 there exists $T>0$ such that
\begin{equation}
\label{eqn:stability-1}
 d^{ss}_{\ch{\Phi}^T(p)}(\ch{\Phi}^T(p),\ch{\Phi}^T(q))
 \leq \frac{1}{4} d^{ss}_p(p,q)
\end{equation}
 for any $p,q \in \ch{M}$.

Let $Y_k$ be the vector field that generates $\ch{\Phi}_k$.
It is H\"older continuous and converges to $\ch{X}$
 as a $C^0$ vector field.
Take $k_0 \geq 1$ so that $Y_k$ is transverse to $\ch{\cF}^{ss}$
 for any $k \geq k_0$.
For any $k \geq k_0$,
 there exists a cocycle $c_k$ over $\ch{\Phi}$ such that
 $\ch{\Phi}_k^{c_k(p,t)}(p) \in \ch{\cF}^{ss}(\ch{\Phi}^t(p))$
 for any $p \in M$ and $t \in \RR$.
We define a flow $\Psi_k$ on $\ch{M}$ by
 $\Psi_k^t(p)=\ch{\Phi}_k^{c_k(p,t)}(p)$.
Since $\ch{\cF}^{ss}$ is H\"older continuous
 and its restriction to each leaf of $\ch{\cF}$ is of class $C^{1+}$,
 we can show that
 the flow $\Psi_k=\{\Psi_k^t\}_{t \in \RR}$
 is of class $C^{1+,0+}_{\ch{\cF}}$
 and the sequence $(\Psi_k)_{k \geq 1}$
 converges to $\ch{\Phi}$ in the $C^{1,0}$-topology.
Hence, there exists $k_1 \geq k_0$ such that
\begin{equation}
\label{eqn:stability-2}
 d^{ss}_{\Psi_k^T(p)}(\Psi_k^T(p),\Psi_k^T(q))
 \leq \frac{1}{2} d^{ss}_p(p,q)
\end{equation}
 for any $p \in \ch{M}$, $q \in \ch{\cF}^{uu}(p)$, and $k \geq k_1$.

Let $\Gamma(\ch{\cF}^{ss})$ be
 the set of continuous self-maps $H$ of $\ch{M}$
 such that $H(p) \in \ch{\cF}^{ss}(p)$ for any $p \in \ch{M}$.
We define a metric $d_\Gamma$ on $\Gamma(\ch{\cF}^{ss})$ by
 $d_\Gamma(H_1,H_2)=\sup_{p \in \ch{M}} d^{ss}_p(H_1(p),H_2(p))$.
For $k \geq k_1$,
 we define a flow $\Gamma_k$ on $\Gamma(\ch{\cF}^{ss})$ by
\begin{equation*}
 \Gamma_k^t(H)(p)=\ch{\Phi}^t \circ H \circ \Psi_k^{-t}(p).
\end{equation*}
The inequality (\ref{eqn:stability-1}) implies
\begin{equation*}
d_\Gamma(\Gamma_k^T(H_1),\Gamma_k^T(H_2))
 \leq \frac{1}{4} d_\Gamma(H_1,H_2).
\end{equation*}
By the contraction mapping principle,
 there exists a fixed point $H_k$ of $\Gamma_k$
 that is homotopic to the identity and satisfies
\begin{equation}
\label{eqn:stability-3}
\sup_{p \in \ch{M}}d^{ss}_p(p,H_k(p)) <C_k
\end{equation}
 for some $C_k>0$.
It is easy to check that $H_k$ is surjective and
 satisfies
\begin{equation*}
 H_k \circ \ch{\Phi}_k^{c_k(p,t)}(p)
 =H_k \circ \Psi_k^t(p)=\ch{\Phi}^t \circ H_k(p)
\end{equation*}
 for any $p \in M$ and $t \in \RR$.
Hence, the proof is completed once we show that $H_k$ is injective.
Suppose $H_k(p)=H_k(q)$ for some $p,q \in \ch{M}$.
Then, both $H_k(p)=H_k(q)$ and $q$ are contained in $\ch{\cF}^{ss}(p)$.
By the inequalities (\ref{eqn:stability-2}) and (\ref{eqn:stability-3}),
 we have
\begin{align*}
 2^l d^{ss}_p(p,q) &\leq
 d^{ss}_{\ch{\Phi}^{-lT}(p)}(\Psi_k^{-lT}(p),\Psi_k^{-lT}(q))\\
 & < d^{ss}_{\ch{\Phi}^{-lT}(p)}
  (H_k \circ \Psi_k^{-lT}(p), H_k \circ \Psi_k^{-lT}(q))
  +2 C_k\\
 &= d^{ss}_{\ch{\Phi}^{-lT}(p)}
  (\ch{\Phi}^{-lT}(H_k(p)),\ch{\Phi}^{-lT}(H_k(q)))+2 C_k
\end{align*}
 for any $l \geq 1$.
By the assumption $H_k(p)=H_k(q)$,
 the last term is equal to $2C_k$,
 and hence, we have $p=q$ by taking $l \ra \infty$.
Therefore, $H_k$ is injective.
\end{proof}

\subsection{Proof of the main theorem}
Fix an $n$-dimensional closed manifold $M$,
 a $C^\infty$ codimension-one
 topologically transitive Anosov flow $\Phi$,
 and a $C^\infty$ metric $g$ on $TM$
 that satisfies (\ref{eqn:Anosov ss,uu}) for some $\lambda>0$.
As the preceding subsection,
 let $\cF$ be the weak stable foliation of $\Phi$
 and $(\nu,\rho_\Phi)$ be the pair
 obtained by applying Theorem \ref{thm:Radon}
 to $\Phi$ and $-(\del \alpha_\Phi/\del t)(\;\cdot\;\cdot\; ;\cF,g)$.

The proof of the main theorem
 is divided into the following two propositions:
\begin{prop}
\label{prop:main 1} 
$\rho_\Phi=1$
 and $\Phi$ is topologically equivalent to a $C^\infty$ Anosov flow
 whose weak unstable foliation is of class $C^{1+}$.
\end{prop}
\begin{prop}
\label{prop:main 2} 
If the weak unstable foliation of $\Phi$ is of class $C^{1+}$,
 then $\Phi$ is topologically equivalent to a $C^\infty$ Anosov
 flow that preserves a $C^\infty$ volume.
\end{prop}
\begin{proof}
[Proof of Propoisition \ref{prop:main 1}]
Let $\cL$ be a $C^{2+}$ one-dimensional foliation
 that is transverse to $\cF$.
As the preceding subsections,
 we take families $\{\vphi_p\}_{p \in M}$
 and $\{\ch{\vphi_p}\}_{p \in M}$ of coordinates
 associated with $(\cF,\cL)$ and $\nu$.
Let $U_p$, $\vphi_{pq}^x$, $\eta_p$, $i_M$, $\ch{M}$, $\ch{\cF}$,
 and $\ch{\Phi}$ be the ones that are given 
 in the preceding subsection.
Remark that the map $\vphi_{pq}^x$ is of class $C^{2+}$.

We say a continuous flow $\Psi$ on a $C^{1+}$ closed
 Riemannian manifold $(M_1,g_1)$
 satisfies {\it the property} $(A)_\epsilon$
 for $\epsilon>0$,
 if there exists a $C^{1+}$ codimension-one foliation $\cF_1$ on $M_1$,
 $T>0$, and $C>0$ such that $\Psi$ is a $C^{1+,0+}_{\cF_1}$ flow and
\begin{align*}
\min\left\{\alpha_{\Psi}^\perp(p,T;\cF_1,g_1),
 -\beta_{\Psi}(p,T;\cF_1,g_1)\right\}
 & > C,\\
 \left|\rho_\Phi \cdot \alpha_{\Psi}(p,T;\cF_1,g_1)
 + \alpha_{\Psi}^\perp(p,T;\cF_1,g_1)\right|
 & < \epsilon \cdot C
\end{align*}
 for any $p \in M_1$.
Since $\alpha_\Psi$ and $\alpha_\Psi^\perp$ are cocycle
 and $\beta_\Psi(p,s+t) \leq
 \beta_\Psi(p,s)+\beta_\Psi(\Psi^s(p),t)$
 for any $s,t \geq 0$,
 we can see that if $\Psi$ satisfies the property $(A)_{\epsilon}$
 and another $C^{1+,0+}$ flow $\Psi'$ on $M_1$
 has the same oriented orbits as $\Psi$,
 then $\Psi'$ also satisfies the property $(A)_\epsilon$.

We will show that $\rho_\Phi=1$
 and that $\Phi$ is topologically equivalent to
 a $C^\infty$ flow on a $C^\infty$ manifold that satisfies
 the property $(A)_{\frac{1}{2}}$.
It completes the proof of the proposition
 by Lemma \ref{lemma:Anosov}.

Fix a finite subset $S$ of $M$ such that $\bigcup_{p \in S}U_p =M$.
Take a $C^{1+}$ partition of unity $\{\zeta_p\}_{p \in S}$
 associated with the covering $\{i_M(U_p)\}_{p \in S}$ of $\ch{M}$.
Let $X$ and $\ch{X}$ be the vector fields
 that generate $\Phi$ and $\ch{\Phi}$ respectively.
For $p \in S$,
 we define a family $(a_{p,m})_{m=1}^{n-1}$ of functions on
 $\ch{\vphi}_p(U_p)$ by
\begin{equation*}
D\ch{\vphi}_p(\ch{X})
 =\sum_{m=1}^{n-1}a_{p,m}\frac{\del}{\del x_m}.
\end{equation*}
We can see
 $D\ch{\vphi}(\ch{X})(x,y)=D\vphi(X)(x,\eta_p^{-1}(y))$
 for any $(x,y) \in \ch{\vphi}_p(U_p)$.
Hence, $a_{p,m}$ and $\del a_{p,m}/\del x_j$
 are H\"older continuous for any $j,m=1,\cdots,n-1$.

Fix a $C^\infty$ non-negative function $\xi$ on $\RR$
 such that its support is compact and its integral over $\RR$ is one.
For $k \geq 1$, we define a function $a_{p,m}^k$ by
\begin{equation*}
a_{p,m}^k(x,y)
 =\int_\RR (\zeta_p \cdot a_{p,m})(x, y+(t/k)) \cdot \xi(t) \;dt.
\end{equation*}
By the standard arguments on mollifiers,
 we can show that $a_{p,m}^k$ is of class $C^{1+}$ and
 it satisfies
\begin{align}
\label{eqn:main 1-1}
\lim_{k \ra \infty}
 \left\|a_{p,m}^k(w) -(\zeta_p \cdot a_{p,m})(w)\right\|_{0,p}
  &=0, \\
\label{eqn:main 1-2}
\lim_{k \ra \infty}
 \sum_{j=1}^{n-1}
\left\|
 \frac{\del a_{p,m}^k}{\del x_j}(w)
 - \frac{\del (\zeta_p \cdot a_{p,m})}{\del x_j}(w)
 \right\|_{0,p}&=0,
\end{align}
 where $\|h\|_{0,p}=\sup_{w \in \ch{\vphi}_p(U_p)}|h(w)|$.

Put
\begin{equation*}
 \ch{X}_k=\sum_{p \in S}
  D\ch\vphi_p^{-1}\left(
   \sum_{m=1}^{n-1}a_{p,m}^k \frac{\del}{\del x_m}
  \right).
\end{equation*}
Then, we have
\begin{align*}
D\ch{\vphi}_q(\ch{X}_k)(x,y)
 &=
 \sum_{p \in S}D(\ch{\vphi}_q \circ \ch{\vphi}_p^{-1})
 \left(\sum_{m=1}^{n-1}
  a_{p,m}^k(x_p,y_p)\frac{\del}{\del x_m}
 \right)\\
 &= \sum_{p \in S}
 \sum_{m=1}^{n-1} a_{p,m}^k(x_p,y_p)
 D(\vphi_{qp}^x)_{x_p}\left(\frac{\del}{\del x_m}\right),
\end{align*}
 where $(x_p,y_p)=\ch{\vphi}_p \circ \ch{\vphi}_q^{-1}(x,y)$.
Since $\vphi_{pq}^x$ is of class $C^{2+}$
 and $\ch{\vphi}_p \circ \ch{\vphi}_q^{-1}$ is of class $C^{1+}$,
 the vector field $D\ch{\vphi}_q(\ch{X}_k)$
 is of class $C^{1+}$ for any $q \in S$.
In particular, $\ch{X}_k$ generates a $C^{1+}$ flow
 on $\ch{M}$ that is tangent to $\ch{\cF}$.
We denote it by $\ch{\Phi}_k$.
By the inequalities (\ref{eqn:main 1-1}) and
  (\ref{eqn:main 1-2}),
 the sequence $(\ch{\Phi}_k)_{k \geq 1}$ converges
 to $\ch{\Phi}$ in the $C^{1,0}$-topology.
It implies that
 the flow $\ch{\Phi}_k$
 is topologically equivalent to $\Phi$ by \ref{lemma:stability}.
For any given $\epsilon>0$,
 the flow $\ch{\Phi}$ satisfies $(A)_\epsilon$.
Hence, $\ch{\Phi}_k$ satisfies the property $(A)_\epsilon$
 if $k$ is sufficiently large.

We show that $\rho_\Phi=1$.
If it is not, there exists $k_0 \geq 1$ such that
 $\ch{\Phi}_{k_0}$
 satisfies the property $(A)_\epsilon$
 for $\epsilon=|1-\rho_\Phi|/2$.
Put $\Psi=\ch{\Phi}_{k_0}$.
Let $E$ be the orthogonal complement of $T\Phi$ 
 with respect to $\ch{g}$. 
Take a metric $\ch{g}'$ on $T\ch{M}$ such that
 the splitting $T\ch{M}=T\Psi \oplus E$
 is also orthogonal with respect to $\ch{g}'$,
 $\ch{g}'|_{E}=\ch{g}|_E$, and $\|X_\Psi\|_{\ch{g}'}=1$,
 where $X_\Psi$ is the vector field that  generates $\Psi$.
Remark that the metric $\ch{g}'$ is H\"older continuous
 and satisfies 
\begin{equation*}
\log{\det}_{\ch{g}'}D\Psi^t(p)
= \alpha_{\Psi}
 (p,t;\ch{\cF},\ch{g})
 +\alpha_{\Psi}^\perp(p,t;\ch{\cF},\ch{g}).
\end{equation*}
Let $T$ and $C$ be constants in the definition 
 of the property $(A)_\epsilon$ for $\Psi$.
Then, we have
\begin{align*}
 \lefteqn{\inf_{p \in \ch{M}}
  \left| \log{\det}_{\ch{g}'}D\Psi^t(p) \right|}\\
 &= \inf_{p \in \ch{M}}\left|\alpha_{\Psi}
 (p,T;\ch{\cF},\ch{g})
 +\alpha_{\Psi}^\perp(p,T;\ch{\cF},\ch{g})
 \right|\\
 & \geq 
 \inf_{p \in \ch{M}}\left\{|1-\rho_\Phi| \cdot
  \left|\alpha_{\Psi}(p,T;\ch{\cF},\ch{g}) \right|
 -\left|\rho_\Phi \cdot \alpha_{\Psi}(p,T;\ch{\cF},\ch{g})
 +\alpha_{\Psi}^\perp(p,T;\ch{\cF},\ch{g}) \right|
 \right\}\\
 & \geq \frac{C}{2}|1-\rho_\Phi|>0.
\end{align*}
It contradicts that $M$ is of finite volume
 with respect to the volume induced from $\ch{g}'$.
So, we have $\rho_Phi=1$.

Fix $k_1 \geq 1$ such that $\ch{\Phi}_{k_1}$ satisfies
 the property $(A)_{(1/2)}$.
By Corollary \ref{cor:foliation},
 $\ch{M}$ admits a $C^\infty$-structure that is compatible
 with the $C^{1+}$-structure and
 such that the orbit foliation of $\ch{\Phi}_{k_1}$
 is generated by a $C^{1+}$ vector field $\ch{Y}$
 as an oriented foliation.
Let $\Psi_1$ be the flow generated by $\ch{Y}$.
It satisfies the property $(A)_{(1/2)}$,

By approximating the vector field $\ch{Y}_1$
 by a $C^\infty$ one,
 we obtain a $C^\infty$ Anosov flow $\Psi_\infty$
 that is topologically equivalent to $\Psi_1$,
 and hence, to $\Phi$.
Since the tangent space of the stable foliation of $\Psi_\infty$
 is $C^0$-close to that of $\Psi_1$,
 we can see that $\Psi_\infty$ satisfies the property $(A)_{(1/2)}$.
As mentioned at the beginning of the proof,
 it completes the proof of the proposition.
\end{proof}

\begin{proof}
[Proof of Proposition \ref{prop:main 2}] 
Let $\cF^u$ be the weak unstable foliation of $\Phi$.
By the assumption,
 we can take the $C^{1+}$ foliation $\cL$ 
 that is transverse to $\cF$ and
 satisfies $\cL(p) \subset \cF^u(p)$ for any $p \in M$.
As the preceding subsection,
 we take families $\{\vphi_p\}_{p \in M}$
 and $\{\ch{\vphi_p}\}_{p \in M}$ of coordinates
 associated with $(\cF,\cL)$ and $\nu$.
Let $\eta_p$, $\ch{M}$, $\ch{\cF}$, and $\ch{\Phi}$
 be the ones that are given in the preceding subsection.

Put $\ch{\cF}^u=i_M(\cF^u)$.
We define a foliation $\cO_p$ on $(-1,1)^{n-1}$ by
\begin{equation*}
\vphi_p^{-1}(\cO_p(x) \times 0)
 \subset \{\Phi^t \circ \vphi_p^{-1}(x,0) \st t \in \RR\}.
\end{equation*}
It is a $C^{1+}$ foliation satisfying
$\vphi_p^{-1}(\cO_p(x) \times (-1,1))
 \subset \cF^u(\vphi_p^{-1}(x,y))$, and hence,
\begin{equation*}
 \ch{\vphi}_p^{-1}(\cO_p(x) \times \eta_p((-1,1)))
 \subset \ch{\cF}^u(\ch{\vphi}_p^{-1}(x,\ch{y})),
 \end{equation*}
 for any $p \in M$, $(x,y) \in (-1,1)^n$,
 and $\ch{y} \in \eta_p((-1,1))$.
Therefore,
 $\ch{\cF}^u$ is a $C^{1+}$ foliation.

Let $\ch{\cO}$ be the orbit foliation of $\ch{\Phi}$.
Each leaf of $\ch{\cO}$ is the intersection of leaves
 of $\ch{\cF}$ and $\ch{\cF}^u$,
 and hence, the foliation $\ch{\cO}$ is of class $C^{1+}$.
By Corollary \ref{cor:foliation},
 $\ch{M}$ admits a $C^\infty$-structure
 that is compatible with the $C^{1+}$-structure of $\ch{M}$
 and such that $\ch{\cO}$ is generated by
 a $C^{1+}$ vector field $Y$.
Let $\Psi$ be the flow generated by $Y$.
Since $\alpha_\Psi$, $\alpha_\Psi^\perp$, and $\beta_\Psi$
 are determined by holonomy maps of $\ch{\cO}$,
 we have
\begin{equation}
\label{eqn:main 2-1}
 \rho_\Phi \cdot \alpha_\Psi^\perp(p,t;\ch{\cF},\ch{g})
 + \alpha_\Psi(p,t;\ch{\cF},\ch{g})=0
\end{equation}
 for any $p \in \ch{M}$ and $t \in \RR$,
 and there exists $T>0$ such that
\begin{equation}
\label{eqn:main 2-2}
\sup_{p \in \ch{M}}\beta_\Psi(p,T;\ch{\cF},\ch{g})<0.
\end{equation}
Since $\alpha_\Psi(p,t) \leq \beta_\Psi(p,t)$,
 Lemma \ref{lemma:Anosov} implies that $\Psi$ is an Anosov flow.

By Proposition \ref{prop:main 1}, we have $\rho_\Phi=1$.
The equations (\ref{eqn:main 2-1}) implies
 that $\Psi$ preserves a H\"older continuous volume.
Now, the proposition follows from Proposition \ref{prop:volume}
 and the structural stability of Anosov flows.
\end{proof}

%
%
\appendix

\section{A short proof of Theorem \ref{thm:foliation}}
\label{sec:A foliation}
In this appendix, we give a short proof of
 Theorem \ref{thm:foliation}.

Let $m$ be the dimension of $\cF$.
We identify the tangent space $T_z \RR^n$ of $\RR^n$
 at $z \in \RR^n$ with $\RR^n=\RR^{m} \oplus \RR^{n-m}$.
Let $\|\cdot\|$ be the Euclidean norm on $\RR^k$ for each $k \geq 1$.
For $\epsilon>0$,
 a $C^\infty$ coordinate $\vphi$ on an open subset $U_\vphi$ of $M$
 is called {\it $(\cF,\epsilon)$-adapted} if
 $\vphi(U_\vphi)=(-9,9)^n$ and
\begin{equation}
\label{eqn:adapted} 
D\vphi(T\cF(p)) \subset
 \{v \oplus v' \in T_{\vphi(p)}\RR^n=\RR^{m} \oplus \RR^{n-m}
 \st \|v'\| \leq \epsilon\|v\| \}
\end{equation}
 for any $p \in U_\vphi$.
For a foliation $\cG$ on $M$,
 let $U^{r+\alpha}(\cG)$ be the set of points $p \in M$
 such that $T\cG$ is of class $C^{r+\alpha}$ on a neighborhood of $p$.
\begin{lemma}
\label{lemma:fol 1}
Let $\vphi$ be an $(\cF,1)$-adapted coordinate.
For any compact subsets $V_0 \subset U^{r+\alpha}(\cF)$
 and any $C^r$-neighborhood $\cU$ of the identity map
 in the space of $C^r$-diffeomorphism of $M$,
 there exists a $C^{r+\alpha}$ diffeomorphism $H \in \cU$ such that
\begin{equation}
\label{eqn:fol 1}
  V_0 \cup \vphi^{-1}([-1,1]^n) \subset U^{r+\alpha}(H(\cF)).
\end{equation}
\end{lemma}
\begin{proof}
By (\ref{eqn:adapted}),
 there exists a $C^{r+\alpha}$ map
 $F:[-4,-4]^n \ra (-9,9)^{n-m}$ such that
 $\vphi^{-1}(w,F(w,z)) \in \cF(\vphi^{-1}(0,z))$
 and $\|F(w,z)-z\| \leq \|w\|$
 for any $(w,z) \in [-4,4]^n =[-4,4]^{m} \times [-4,4]^{n-m}$.
Put $U =\vphi^{-1}(\{(w,F(w,z)) \st (w,z) \in [-4,4]^n\})$.
For a map $G$ from $[-4,4]^n$ to $\RR^n$ or $\RR^{n-m}$,
 let $V^{r+\alpha}(G)$ be the set of points $w \in \RR^n$
 such that $G$ is of class $C^{r+\alpha}$ on a neighborhood of $w$.

Fix a $C^\infty$ function $\lambda$ on $[-4,4]^n$ such that
 $0 \leq \lambda(w,z) \leq 1$ for any $(w,z) \in [-4,4]^n$,
 $\lambda(w,z)=0$ on $[-4,4]^n \-- [-7/2,7/2]^n$,
 and $\lambda(w,z)=1$ on $[-3,3]^n$.
For any given $C^\infty$ map $f:[-4,4]^n \ra (-9,9)^{n-m}$,
 we define $C^{r+\alpha}$ maps $G_f:[-4,4]^n \ra (-9,9)^n$
 and $H_f:M \ra M$ by
\begin{equation*}
 G_f(w,z)=(w,(1-\lambda(w,z))F(w,z)+\lambda(w,z)f(w,z))
\end{equation*}
 and
\begin{equation}
\label{eqn:fol 1-0}
 H_f(p) =\left\{
 \begin{array}{ll}
 \vphi^{-1} \circ G_f \circ G_F^{-1} \circ \vphi(p) & (p \in U)\\
 p & (p \not\in U).
 \end{array}
 \right.
\end{equation}
Remark that $[-1,1]^n \subset G_F((-3,3)^n)$ and
\begin{equation}
\label{eqn:fol 1-1}
 V^{r+\alpha}\left(\frac{\del F}{\del x_j}\right) \cup [-3,3]^n
  \subset V^{r+\alpha}\left(\frac{\del G_f}{\del x_j}\right)
\end{equation}
 for any $j=1,\cdots,n$.

If $f$ is sufficiently $C^1$-close to $F$,
 then $H_f$ is a diffeomorphism.
Since
\begin{displaymath}
 \vphi^{-1} \circ G_f([-4,4]^{m} \times z) \subset
 H_f(\cF)(\vphi^{-1} \circ G_f(w,z))
\end{displaymath}
 for any $(w,z) \in [-4,4]^n$,
 we have
\begin{equation}
\label{eqn:fol 1-2}
 U^{r+\alpha}(H_f(\cF)) \cap H_f(U)
  = \vphi^{-1} \circ G_f\left(\bigcap_{j=1}^n V^{r+\alpha}
  \left(\frac{\del G_f}{\del x_j}\right)\right).
\end{equation}
In particular,
\begin{equation}
\label{eqn:fol 1-3}
 U^{r+\alpha}(\cF) \cap U
  = \vphi^{-1} \circ G_F\left(\bigcap_{j=1}^n V^{r+\alpha}
  \left(\frac{\del F}{\del x_j}\right)\right)
\end{equation}
 since $G_F(x,y)=(x,F(x,y))$ and $H_F$ is the identity map.

There exists a $C^r$-neighborhood ${\cal V}$ of $F$
 such that $H_f$ is a diffeomorphism in $\cU$,
 $V_0 \subset H_f(U^{r+\alpha}(\cF))$,
 and $[-1,1]^n \subset G_f((-3,3)^2)$ for any $f \in {\cal V}$.
By the relations (\ref{eqn:fol 1-0}),
 (\ref{eqn:fol 1-1}),  (\ref{eqn:fol 1-2}),
 and  (\ref{eqn:fol 1-3}),
 the diffeomorphism $H_f$ satisfies
 the required condition (\ref{eqn:fol 1})
 for any $C^\infty$ map $f \in {\cal V}$.
\end{proof}

Now, we prove Theorem \ref{thm:foliation}.
Fix a family $\{\vphi_i\}_{i=1}^k$
 of $(\cF,1/2)$-adapted coordinates so that
 $\bigcup_{i=1}^k \vphi_i^{-1}([-1,1]^n)=M$.
Take a $C^r$ neighborhood $\cU$ of the identity map on $M$
 so that $\vphi_i$ is $(H(\cF),1)$-adapted
 for any $H \in \cU$ and $i=1,\cdots,k$.

Proof is by induction of $i$.
By Lemma \ref{lemma:fol 1},
 there exists $H_1 \in \cU$ such that
 $\vphi_1^{-1}([-1,1]^n) \subset U^{r+\alpha}(H_1(\cF))$.
Suppose that there exists $H_i \in \cU$ such that
 $\bigcup_{j=1}^i \vphi_j^{-1}([-1,1]^n) \subset
 U^{r+\alpha}(H_i(\cF))$ for some $i \leq k-1$.
Applying Lemma \ref{lemma:fol 1} to $\vphi_{i+1}$ and
 $\bigcup_{j=1}^i \vphi_i^{-1}([-1,1]^n)$,
 we obtain a diffeomorphism $H_{i+1} \in \cU$
 such that $\bigcup_{j=1}^{i+1} \vphi_j^{-1}([-1,1]^n) \subset
 U^{r+\alpha}(H_{i+1}(\cF))$.
It completes the induction.

%
%
%


\begin{thebibliography}{99}
\bibitem{Ar}

A.Arbieto, C.Matheus,
 A pasting lemma and some applications for conservative systems.
 {\it Ergodic Theory Dynam. Systems} {\bf 27} (2007), no. 5, 1399--1417.

\bibitem{Ba}
T.Barbot, G\'eom\'etrie transverse des flots d'Anosov.
 Th\'ese, Universit\'e de Lyon I, 1992.


\bibitem{Bo}
R.Bowen,
Equilibrium states and the ergodic theory of Anosov diffeomorphisms. 
 Lecture Notes in Mathematics, {\bf 470}.
 Springer-Verlag, Berlin-New York, 1975.

\bibitem{Ca}
E.~E.~Cawley,
The Teichm\"uller space of an Anosov diffeomorphism of $T^2$.
 {\it Invent. Math.} {\bf 112} (1993), no. 2, 351--376.


\bibitem{Do}
C.I.Doeling,
 Persistently transitive vector fields on three-dimensional manifolds.
 {\it Dynamical systems and bifurcation theory (Rio de Janeiro, 1985)},
 59--89, Pitman Res. Notes Math. Ser., {\bf 160},
 {\it Longman Sci. Tech., Harlow}, 1987.

\bibitem{FW}
J.Franks, R.Williams,
 Anomalous Anosov flows.
 {\it Global theory of dynamical systems
 (Proc. Internat. Conf., Northwestern Univ., Evanston, Ill., 1979)},
 pp. 158--174, 
 Lecture Notes in Math., {\it 819}, Springer, Berlin, 1980.

\bibitem{Gh}
E.Ghys,
 Codimension one Anosov flows and suspensions.
 {\it Dynamical systems, Valparaiso 1986},
 59--72, Lecture Notes in Math., {\it 1331}, Springer, Berlin, 1988.

\bibitem{Ha}
D.Hart,
On the smoothness of generators.
 {\it Topology} {\bf 22} (1983), no. 3, 357--363.

\bibitem{KH}
 A.~Katok and B.~Hasselblatt,
 Introduction to the modern theory of dynamical systems.
 Encyclopedia of Mathematics and its Applications,
 54. {\it Cambridge University Press, Cambridge}, 1995.

\bibitem{Li}
Liv\v sic, A. N.,
Certain properties of the homology of $Y$-systems. 
{\it Math. Notes} {\bf 10} (1971), 758--763.
(translated from {\it Mat. Zametki} {\bf 10} (1971), 555--564)

\bibitem{LS}
A.N.Liv\v sic and Ja.G.{Sina\u \i},
Invariant measures that are compatible
 with smoothness for transitive $C$-systems. (Russian) 
{\it Dokl. Akad. Nauk SSSR} {\bf 207} (1972), 1039--1041.

\bibitem{LMM}
R.~de la Llave, J.M.~Marco,
 and R.Moriy\'on,
 Canonical perturbation theory of Anosov systems
 and regularity results for the Liv\v sic cohomology equation. 
 {\it Ann. of Math. (2)} {\bf 123} (1986), no. 3, 537--611.

\bibitem{Ra}
M.Ratner,
Markov partitions for Anosov flows on $n$-dimensional manifolds. 
{\it Israel J. Math.} {\bf 15} (1973), 92--114.

\bibitem{Sh}
  M.~Shub,
  {\it{Global stability of dynamical systems}},
  {\it Springer-Verlag, Berlin-New York}, 1986.

\bibitem{Si}
 S.N.Simi\'c,
 Volume preserving codimension one Anosov flows
 in dimensions greater than three are suspensions.
 preprint, arXiv:math.DS/058024, 2005.

\bibitem{Ve}
A.Verjovsky,
Codimension one Anosov flows. 
{\it Bol. Soc. Mat. Mexicana (2)} {\bf 19} (1974), no. 2, 49--77.

\end{thebibliography}
\end{document}